\documentclass[12pt]{amsart}
\usepackage{amssymb}
\usepackage{amsmath}
\usepackage{hyperref}
\newcommand{\R}{\mathbb{R}}

\font\eufm=eufm10
\def\frak#1{\hbox{\eufm#1}}
\newcommand{\bd}{\begin{document}}
\newcommand{\ed}{\end{document}}
\newcommand{\be}{\begin{enumerate}}
\newcommand{\ee}{\end{enumerate}}
\newcommand{\bi}{\begin{itemize}}
\newcommand{\ei}{\end{itemize}}
\newcommand{\ba}{\begin{array}}
\newcommand{\ea}{\end{array}}
\newcommand{\vs}{\vspace*{0.3\baselineskip}}
\newcommand{\vsm}{\vspace*{-0.3\baselineskip}}
\newcommand{\kom}[1]{{\bf #1}}
\newcommand{\komm}[1]{}
\newtheorem{defi}{Definition}[section]
\newtheorem{tw}[defi]{Theorem}
\newtheorem{prop}[defi]{Proposition}
\newtheorem{lem}[defi]{Lemma}
\newtheorem{re}[defi]{Remark}
\newtheorem{col}[defi]{Corollary}
\newtheorem{ex}[defi]{Examples}
\newtheorem{ex1}[defi]{Example}
\newtheorem{zad}{Exercise}[section]
\newtheorem{zal}{Assumptions}[section]
\newtheorem{assumpt}[defi]{Assumptions}
\newcommand{\Om}{\Omega}
\newcommand{\om}{\omega}
\newcommand{\G}{\Gamma}
\newcommand{\D}{\Delta}
\renewcommand{\d}{\delta}
\newcommand{\ga}{\gamma}
\newcommand{\eps}{\epsilon}
\newcommand{\ove}{\overline}
\newcommand{\ms}{\oplus}
\newcommand{\mt}{\otimes}
\newcommand{\dz}{\wedge}
\newcommand{\lra}{\longrightarrow}
\newcommand{\sign}{\mbox{$ sgn $}}
\newcommand{\rel}{\mbox{$\,$\rule[0.5ex]{1.1em}{0.2pt}$\triangleright\,$}}
\newcommand{\dow}{\hspace*{\fill}\rule{1.6ex}{1.6ex}\hspace*{0.1em}}
\newcommand{\dowl}{\hspace*{\fill}\rule{1ex}{1ex}\hspace*{1em}}
\newcommand{\sd}{\hspace{0.3ex}\tiny{\rhd\mbox{\hspace{-2ex}}<}\hspace{0.3ex}}
\newcommand{\mmt}[2]{\mbox{$\vphantom{}_{#1}\times_{#2}$}}
\newcommand{\gotg}{\frak g}
\newcommand{\gota}{\frak a}
\newcommand{\gotb}{\frak b}
\newcommand{\gotc}{\frak c}
\newcommand{\gothh}{\frak h}
\newcommand{\gott}{\frak t}
\newcommand{\hd}{\hat{\d}}
\newcommand{\oml}{\Omega_L^{1/2}}
\newcommand{\omr}{\Omega_R^{1/2}}
\newcommand{\omh}{\Omega^{1/2}}
\newcommand{\lo}{\lambda_0}
\newcommand{\ro}{\rho_0}
\newcommand{\lma}{\Lambda^{max}}
\newcommand{\timh}{\times_h}
\newcommand{\Gd}{\G^{(2)}}
\newcommand{\el}{e_L}
\newcommand{\er}{e_R}
\newcommand{\GG}{\G_1\times\G_2}
\newcommand{\gdot}{\hspace{-0.1em}\cdot\hspace{-0.1em}}
\newcommand{\tran}{\frown\hspace{-2.2ex}|\hspace{1.9ex}}
\def\tgr{{\bf t}}
\def\sgr{{\bf s}}
\def\fgr{{\bf f}}
\def\rhogr{{\boldsymbol \rho}}
\newcommand{\la}[2]{\Lambda_{#1#2}}
\newcommand{\ad}{\mathrm{ad}}
\newcommand{\kad}{\ad^{\#}}
\newcommand{\wl}[1]{\vphantom{X}_{#1}{\G}}
\newcommand{\te}{\tilde{e}}
\newcommand{\notka}[1]{}
\newcommand{\sA}{\mbox{$\mathcal A$}}
\newcommand{\sT}{\mbox{$\mathcal T$}}
\newcommand{\sB}{\mbox{$\mathcal B$}}
\newcommand{\sE}{\mbox{$\mathcal E$}}
\newcommand{\sF}{\mbox{$\mathcal F$}}
\newcommand{\sL}{\mbox{$\mathcal L$}}
\newcommand{\sO}{\mbox{$\mathcal O$}}
\newcommand{\sD}{\mbox{$\mathcal D$}}
\newcommand{\sM}{\mbox{$\mathcal M$}}
\newcommand{\sS}{\mbox{$\mathcal S$}}
\newcommand{\sY}{\mbox{$\mathcal Y$}}
\newcommand{\hY}{\mbox{$\hat{Y}$}}
\newcommand{\hS}{\mbox{$\hat{S}$}}
\newcommand{\hX}{\mbox{$\hat{X}$}}
\newcommand{\dif}{differential }
\newcommand{\gru}{groupoid }
\newcommand{\grus}{groupoids }
\newcommand{\ti}{\tilde}
\newcommand{\halden}{half density }
\newcommand{\haldens}{half densities }
\renewcommand{\top}{topological }
\newcommand{\Setrel}{\mbox{\rm SetRel}}

\newcommand{\cstardwa}{\mbox{$C^*_r(\Gamma\times\Gamma)$}}
\newcommand{\skal}{\,|\,}
\newcommand{\hil}[1]{\mbox{$\mathcal{#1}$}}
\newcommand{\cc}{C_{fin}}
\newcommand{\rra}{\rightrightarrows}
\newcommand{\Hh}{H_0}
\newcommand{\sect}{{\mathrm {Sec}}}
\newcommand{\LL}{L^2}
\newcommand{\skalp}[2]{#1\hspace{-0.4ex} \cdot \hspace{-0.4ex} #2}
\newcommand{\mskalp}[2]{\mbox{$\displaystyle \left<\, #1\skal #2 \,\right>$}}
\newcommand{\wtil}{\widetilde}
\newcommand{\dsp}{\displaystyle}
\newcommand{\npj}{n\hspace{-0.5ex}+\hspace{-0.5ex}1}
\newcommand{\npji}{n\hspace{-0.3ex}+\hspace{-0.3ex}1}
    \newcommand\lie[1]{\sL(#1)}
    \newcommand\bunR{\widetilde{\R}}
\begin{document}
\title{A quantum space of Euclidean lines.}
\author{Piotr Stachura}
\address{Institute of Information Technology, Warsaw University of Life Sciences -  SGGW\\Warszawa, Poland.
  \thanks{e-mail:piotr\_stachura1@sggw.edu.pl}}
\date{}
\begin{abstract} This article presents a differential groupoid with ``coaction'' of the groupoid underlying the
  Quantum Euclidean Group (i.e. its $C^*$-algebra is the  $C^*$-algebra of this quantum group).
  The dual of the Lie algebroid is a Poisson manifold that can be identified with the space of oriented
  lines in Euclidean space equipped with a Poisson action of the Poisson-Lie Euclidean group.
  \end{abstract}
\maketitle
\section{Introduction}

The inspiration for  this work was  the article of A. Ballesteros, I. Gutierrez-Sagredo and F.J. Herranz \cite{BalGH}.
They found a Poisson structure on the space of timelike, oriented, future (or past) pointing straight lines in Minkowski space by considering it as a
homogeneous space of the Poincar\'{e} Group.   They observed also that stabilizer of a
timelike line is a Poisson-Lie subgroup of a Poisson-Lie Poincar\'{e} Group related to the $\kappa$-deformation.

This suggests the similar result for $C^*$-version i.e. $\kappa$-Poincar\'{e} \cite{PS-kappa2},
namely that this stabilizer is a quantum subgroup and  that there is the  corresponding quantum homogeneous space.
It seems everything goes smoothly  at least for quantum groups  defined by Double Lie Groups \cite{PS-DLG}
(special case of bicrossed product construction of \cite{Va-Ve}) and will be presented in the separate article \cite{PS-future}.
\komm{Przypomniec sobie te koizotropowe i poissona -liego podalgebry}

The example presented here is related to the Quantum Euclidean Group and  should be considered as ``the quantum space of (oriented) euclidean lines''.
The euclidean case is technically simpler but  similar construction works for $\kappa$-Poincar\'{e}; it will be presented in \cite{PS-future}
as well as a full justification of its name as a quantum homogeneous space for a  locally compact quantum group as described in \cite{DKSS, Kasp-Sol}.
Here only a partial justification for the name is provided -- the semi-classical limit.

The $\npj\,(n\geq 2)$ dimensional Quantum Euclidean Group (in the following abbreviated as $QE(\npj)$) is defined by the Iwasawa decomposition of
$G:=SO_0(1,\npj)$ and its $C^*$-algebra is $C^*(G_B)$ -- the groupoid $C^*$-algebra of $G_B:G\rra SO(\npj)$ (see sections  \ref{sect:DG}, \ref{sect:QE} for notation). 
\komm{UWAGA: Role $\tilde{Z}$ oraz $Z$ i $\tilde{T}$ oraz $T$ ZAMIENIONE stosunku do subgroups!!!!} 
It turns out there exists structure of a differential groupoid $Z:TS^n\times\R_+\rra S^n$ and a (Zakrzewski) morphism, ``coaction'' of $G_B$ on $Z$, 
$\delta_{Z}:Z\rel G_B\times Z$ satisfying $(\delta_B\times id)\delta_{Z}=(\delta_{Z}\times id)\delta_{Z}$,
where $\delta_B:G_B\rel G_B\times G_B$ is the  coassociative morphism whose  $C^*$-lift defines  the   comultiplication of $QE(\npj)$.
The submanifold  $TS^n\simeq TS^n\times \{1\}\subset TS^n\times\R_+$ is a (wide) subgroupoid $\tilde{S}:TS^n\rra S^n$ of $Z$
and $\delta_{Z}$ restricts to $\tilde{S}$, and defines the  morphism $\delta_{\tilde{S}}:\tilde{S}\rel G_B\times \tilde{S}$.
Applying cotangent lifts \cite{SZ2} to
morphisms $\delta_{Z}$ and  $\delta_{\tilde{S}}$ one obtains morphisms of symplectic groupoids
$T^*\delta_{Z}: T^*Z\rel T^*G_B\times T^*Z$, and $T^*\delta_{\tilde{S}}: T^*\tilde{S}\rel T^*G_B\times T^*\tilde{S}$.
Base maps of these morphisms are Poisson maps and actions of the Poisson-Lie Euclidean Group on $\lie{Z}^*$ and $\lie{\tilde{S}}^*$ respectively,
where e.g. $\lie{Z}^*$ is the dual bundle of the Lie algebroid of $Z$ (so a Poisson manifold in the canonical way). 

The left(target)  projection in $Z$ is just the bundle projection in $TS^n$ therefore its Lie algebroid $\lie{Z}$
is just the bundle $TS^n\oplus \bunR$, where $\bunR$ is the trivial bundle $S^n\times \R$, and $\lie{\tilde{S}}$ is  $TS^n$.
As a vector bundle $\lie{Z}=TS^n\oplus \bunR$ is trivial; it can be identified with $S^n\times \R^{n+1}$ 
(the bundle  $\bunR=S^n\times \R$ is identified with the normal bundle to $S^n$ in $\R^{n+1}$). By using the euclidean structure, bundles
$TS^n$ and $T^*S^n$ can be identified, this way one obtains Poisson structures on $TS^n$ and $TS^n\oplus\bunR=S^n\times \R^{\npji}$ with Poisson
actions of (Poisson-Lie) Euclidean Group.

On the other hand for the affine Euclidean Space $(M,V)$ of dimension $\npj$,
the Euclidean Group act transitively on the sphere bundle over $M$: $S^n\times M$ and on the space of
oriented straight lines, so both can be identified with homogeneous spaces after a choice of a point $(p_0,m_0)\in S^n\times M$ or a line determined by
$(p_0, m_0)$. It turns out, that after these identifications (and after the choice of orthonormal basis in $V$) the sphere bundle is
$S^n\times\R^{\npji}$ and the space of lines is  $TS^n$ and these two actions coincide with actions obtained from groupoids.

Let me make two comments on the article \cite{BalGH}. The first one is a small correction: what is considered there, strictly speaking, is not the space of
``oriented, timelike lines'' but only ``half of it'' i.e. future or past pointing. This is clear from the construction:
dividing  connected component of the identity of Poincar\'{e} Group by stabilizer of the line,  and this connected component does not act
transitively on the set of oriented  timelike lines. And the second is that 
Authors of  \cite{BalGH} observed their  Poisson bracket is symplectic off some three dimensional submanifold (they write off the point, but it seems 
they mean a point in ``half of coordinates''), and provide symplectic coordinates there.
The similar situation is in our construction --  the groupoid $\tilde{S}$ has two orbits: the distinguished point $p_0$ and
$S^n\setminus\{p_0\}$;  the transitive component over $S^n\setminus\{p_0\}$ has no isotropy so it is (isomorphic to) the {\em pair groupoid}
$\R^n\times \R^n\rra \R^n$, and  it is clear that the corresponding Poisson bracket is symplectic there. 

It seems they ignore the fact that their Poisson structure is not symplectic everywhere and they write in the introduction:
{\em As the most relevant result, we prove in Section 5 that the non-trivial Poisson structure on the space of worldlines so obtained
  is just a canonical symplectic structure} and in section 5 (after formula (32)):
{\em since what we have found is that the homogeneous Poisson structure induced by the $\kappa$-Poincar\'{e}
  r-matrix on the space or worldlines is just a symplectic structure on W (without the origin)}.
But they {\em have not shown}  their structure is symplectic everywhere and most probably {\em it is not}, and this {\em makes a difference } for
the structure of a differential groupoid $C^*$-algebra that should be regarded as quantization the  Poisson structure considered.

Let us briefly describe the content of each section.
The next one  establishes  notation and basic facts about groupoids related to group decompositions.
It also contains construction of division of a groupoid by a group of automorphism with groupoid structure on the set of equivalence classes.
To my  best knowledge this was not known before although after that discovery I found related construction in \cite{Tayl}.
There is a remark, how this two are related.
In the third, short section I recall classical description of the space of oriented lines in Euclidean space as a homogeneous space of Euclidean Group.
In the fourth one the main construction is applied to the Quantum Euclidean Group and   its semi-classical limit is computed.

\section{Groupoids related to group decomposition --  algebra.}\label{sect:DG}
\subsection{Groupoids} The Zakrzewski Category of groupoids will be used in the following (see \cite{SZ1} for short or  \cite{PS-short} for more detailed  exposition).
It is defined using relations, so let us recall basic facts and establish needed notation. Symbol $r: X\rel Y$ stands for {\em a relation}
from $X$ to $Y$; usually relation will be identified with its graph i.e. a subset of $Y\times X$ -- notice the order of factors;
{\em a domain } of $h$ is denoted by $D(h)$. If $r:X\rel Y$ is a relation, $r^T:Y\rel X$ denotes its {\em transposition} : $(y,x)\in r\iff(x,y)\in r^T$.
Transposition is an involution and reverses order of composition: $(r^T)^T=r\,,\,(r s)^T=s^T r^T$.

A groupoid $\Gamma$ with the set of unit $E\subset \Gamma$ will be denoted by $\Gamma\rightrightarrows E$.
The source and target projections will be denoted by $\er$ and $\el$ resp.; a groupoid multiplication relation by $m:\Gamma\times\Gamma\rel \Gamma$, and
a groupoid inverse will be denoted by $s$; if structure relations should be explicitly mentioned a groupoid will be denoted by $(\Gamma, m,s,E)$.
For $\gamma\in\Gamma$ we will use notation $F_l(\gamma),F_r(\gamma)$ for left and right
fiber passing through $\gamma$, i.e.  $F_l(\gamma):=\el^{-1}(\el(\gamma))$.
In Zakrzewski Category of groupoids   morphism from
$\Gamma_1$ to $\Gamma_2$ is a {\em relation} $h:\Gamma_1\rel\Gamma_2$ that preserves groupoid structure:
$hm_1=m_2(h\times h)\,,\, hs_1=s_2h,\,,\, h(E_1)=E_2$.  A subset $\tilde{\Gamma}\subset \Gamma$  is a {\em subgroupoid} iff
$m(\tilde{\Gamma}\times \tilde{\Gamma})\subset \tilde{\Gamma}$ and $s(\tilde{\Gamma})=\tilde{\Gamma}$;  a subgroupoid $\tilde{\Gamma}\subset \Gamma$
is {\em wide} iff $E\subset \tilde{\Gamma}$.
The flip map $(x,y)\mapsto (y,x)$ is  denoted \mbox{by $\sim$.}
\subsection{Groupoids related to group decomposition} It is useful to introduce the following:
\notka{def:tis}

\begin{defi}\label{def:tis}
 {\em  Let $G$ be a group and  $A, C\subset G$ be subgroups with $A\cap C=\{e\}$.
   In this situation the triple $(G;A,C)$ will be called } a trivially intersecting subgroups {\em and abbreviated as TIS. If, additionally, $AC=G$ the
   triple $(G; A,C)$ is called} a double group {\em and abbreviated as DG}. {\em A double group  $(G; A,C)$ with $G$ a Lie group and $A,C$ closed subgroups
 will be called} a double Lie group \cite{LuWe}  {\em and abbreviated as DLG}.
  \end{defi}

Let $(G; A, C)$ be  TIS.  Any element $g$ in the set $\G:=AC\cap CA$ of {\em decomposable elements} has  unique decompositions
  $g=a_L(g)c_R(g)=c_L(g) a_R(g)\,,\,a_L(g),a_R(g)\in A\,,\,  c_R(g), c_L(g)\in C$. The set $\G$ carries two groupoid structures
  denoted by $\Gamma_A:\Gamma\rra A$ and $\Gamma_C:\Gamma\rra C$. The inverse $s_A$ and multiplication relation
  $m_A:\Gamma_A\times \Gamma_A\rel\Gamma_A$ are given by:
\begin{align*} s_A(g)&:= c_L(g)^{-1}a_L(g)=a_R(g) c_R(g)^{-1}\,,\,\,\quad
Gr(m_A) := \{(c_1 a c_2; c_1 a, a c_2) : c_1a, ac _2\in\Gamma\}
\end{align*}
and similarly for $\Gamma_C$.  If $(G; A, C)$ is DG (i.e. $\Gamma=G$) these groupoids will be denoted by $G_A$ and $G_C$ (see \cite{PS-kappa2} for details).
The relation $\delta_A:=m_C^T: \Gamma_A\rel \Gamma_A\times\Gamma_A$ is a
{\em coassociative relation} i.e. $(id\times \delta_A)\delta_A=(\delta_A\times id)\delta_A$,  but {\em not a morphism of groupoids}.
It satisfies equalities:\notka{eq:amor-delta}
\begin{equation} \label{eq:amor-delta} \delta_A(A)=A\times A\,,\quad \delta_A s_A=(s_A\times s_A) \delta_A
  \end{equation}
but instead of the third equality relating multiplications there is only the inclusion:\notka{eq:amor-delta1}
\begin{equation}\label{eq:amor-delta1}
  \delta_A m_A\supset (m_A\times m_A)(id \times \sim \times id) (\delta_A\times \delta_A),\quad {\rm where }\quad \sim: (x,y)\mapsto (y,x)
  \end{equation}
and the equality is equivalent to  $AC=CA$ i.e. to $(G;A,C)$ being DG (see \cite{SZ1}, Prop.9.2).

Let $(G;A,C)$ be a TIS, $N(C):=\{g\in G: gCg^{-1}=C\}$ be the {\em normalizer} of $C$ and  
let $A_0:=N(C)\cap A$;  clearly $A_0$ is a subgroup and $CA_0=A_0C\subset \Gamma $, and moreover we have: \notka{prop:automor-A0}
\begin{prop} \label{prop:automor-A0} Let $(G; A, C)$ be TIS and  $A_0:=A\cap N(C)$. Then
  \begin{enumerate}
    \item $A_0\Gamma=\Gamma A_0=\Gamma$ i.e. multiplication by $A_0$ preserves the set of decomposable elements.
    \item  For  $a\in A_0$ the mappings
  $$L_a: \Gamma_A\ni\gamma \mapsto a \gamma\in\Gamma_A\,,\,\, R_a: \Gamma_A\ni\gamma \mapsto  \gamma a \in\Gamma_A $$
  are groupoid automorphisms.
  \item For $a\in A_0$:  $\displaystyle \delta_A R_a=(id\times R_a)\delta_A$ 
 \end{enumerate}
 \end{prop}
 \noindent {\em Proof: } (1) Since $A_0$ is a group it is sufficient to show that $A_0\Gamma\subset \Gamma$ and $\Gamma A_0\subset \Gamma$.
 Let  $\gamma=a c=\tilde{c} \tilde{a}$ and $a_0\in A_0$. Then 
 $\gamma a_0= \tilde{c} (\tilde{a} a_0) \in CA$ and  $\gamma a_0= aa_0 (a_0^{-1} c a_0) \in A C$, and similarly for $a_0\gamma$.\\
(2) Since, for $a\in A_0$: $(L_a)^{-1}=L_{a^{-1}}$ and $(R_a)^{-1}=R_{a^{-1}}$ it is sufficient to show that $L_a, R_a$ are
morphisms $\Gamma_A\rightarrow \Gamma_A$ i.e. $L_a(A)=A$ and $L_{a}s_A=s_AL_{a}$,  and  $m_A (L_{a}\times L_{a})=L_{a} m_A$, and similarly for $R_a$.

a)  Clearly for any $a\in A_0$ one has  $L_a(A)=R_a(A)=A$;

b)  Let $\gamma=ac=\tilde{c}\tilde{a}$ and $a_0\in A_0$, then clearly $a_L(a_0\gamma)=a_0 a_L(\gamma)$ and
$a_R(\gamma a_0)= a_R(\gamma) a_0$, and 
$$a_R(a_0\gamma)=a_R((a_0\tilde{c}a_0^{-1}) a_0 \tilde{a})=a_0 a_R(\gamma)\,,\,a_L(\gamma a_0)=a_L(aa_0(a_0^{-1}ca_0))=a  a_0 =a_L(\gamma)a_0$$
Therefore
$$a_0 s_A(\gamma)=a_0 a_R(\gamma) c_R(\gamma)^{-1}=a_R(a_0\gamma) c_R(a_0 \gamma)^{-1}=s_A(a_0\gamma),$$
$$s_A(\gamma)a_0= c_L(\gamma)^{-1} a_L(\gamma) a_0=c_L(\gamma a_0)^{-1}  a_L(\gamma a_0)=s_A(\gamma a_0),$$
i.e. $L_{a_0}s_A=s_AL_{a_0}$ and  $R_{a_0}s_A=s_A R_{a_0}$;

c) For any $\gamma_1, \gamma_2, \gamma_3\in\Gamma_A$ :
$$(\gamma_1;\gamma_2, \gamma_3)\in L_{a_0} m_A \iff a_R(\gamma_2)=a_L(\gamma_3) {\rm \,and\,} \gamma_1=a_0 \gamma_2 c_R(\gamma_3)$$
\begin{equation*}\begin{split}
    (\gamma_1;\gamma_2, \gamma_3)\in m_A (L_{a_0}\times L_{a_0}) \iff &  a_R(a_0 \gamma_2)=a_L(a_0 \gamma_3) {\rm \,and\,} \gamma_1=a_0\gamma_2 c_R(a_0\gamma_3)
    \\ \iff&  a_R(\gamma_2)=a_L(\gamma_3) {\rm \,and\,} \gamma_1=a_0\gamma_2 c_R(\gamma_3), \end{split}\end{equation*}
where for the second equivalence b) was used. The corresponding equality for $R_{a_0}$ may be proven in the same way.

\noindent (3) Let $a\in A_0$, by definition of $\delta_A$:
$$(\gamma_1,\gamma_2;\gamma_3)\in \delta_A R_a\,\iff\, c_R(\gamma_1)=c_L(\gamma_2) \, {\rm and\,} \gamma_3 a=\gamma_1 a_R(\gamma_2), $$
$$(\gamma_1,\gamma_2;\gamma_3)\in (id\times R_a)\delta_A \,\iff\, c_R(\gamma_1)=c_L(\gamma_2a^{-1}) \, {\rm and\,} \gamma_3 =\gamma_1 a_R(\gamma_2a^{-1}),$$
and these  two are equivalent because of (1) and $c_L(\gamma_2a^{-1})=c_L(\gamma_2)\,,\,a_R(\gamma_2a^{-1})=a_R(\gamma_2) a^{-1}$.

\dowl

 Let us now recall  some properties of equivalence relations
  which will be used in the following. Let $\pi:X\rightarrow Y$ be a surjection,  then $h:=\pi^T\pi:X\rel X$ is the
  equivalence relation on $X$ defined by $\pi$: $(x,x')\in h \iff \pi(x)=\pi(x')$; moreover  
  $\pi\pi^T=id_Y$  and $h^T h=\pi^T\pi \pi^T\pi=\pi^T\pi=h$. Since any equivalence relation is of this form
  ($\pi$ is the canonical surjection on the set of equivalence classes) one obtains the characterization: \notka{eq:rel-row}
   \begin{equation}\label{eq:rel-row}
     h:X\rel X  \,\,\mbox{\rm is an equivalence relation} \,\iff\, D(h)=X\,\,{\rm and}\,\, h=h^Th
   \end{equation}
   Notice also that for any $h:X\rel X$ the equality $h=h^Th$ is equivalent to $h=h^T=h^2$, and (\ref{eq:rel-row}) implies
   that if $h_i:X_i\rel X_i,\, i=1,2$ are equivalence relations so is $h_1\times h_2:X_1\times X_2\rel X_1\times X_2$.
   Writing equivalence relation $h:X\rel X$ as $h=\pi^T\pi$ 
   one easily proves:\notka{lem:hkomut}
   \begin{lem}\label{lem:hkomut}
     Let $h_X:X\rel X$ and $h_Y:Y\rel Y$ be equivalence relations and $\pi_X, \pi_Y$ canonical projections. Assume a relation $k:X\rel Y$ satisfies
     $h_Y k\subset k h_X$. Then $h_Y k \pi_X^T =\pi_Y^T \pi_Y k \pi_X^T=k \pi_X^T$.
   \end{lem} \dowl
   
\noindent And for group actions we have: 
\begin{lem}
  Let a group $A$ act on sets $X$ and $Y$ (from left or right) by $\varphi_a:X\rightarrow X\,,\,\psi_a:Y\rightarrow Y\,,\,a\in A$.
  Let $k:X\rel Y$ be a relation satisfying $k \varphi_a=\psi_a k$ for every  $a\in A$ and 
  let $h_X$ and $h_Y$ be corresponding  equivalence relations.  
  Then $k h_X=h_Y k.$
\end{lem}
\dowl

Let $Z:=\Gamma\slash A_0$ (via the right action  $R$, as in prop. \ref{prop:automor-A0})  and let $\pi_Z:\Gamma\rightarrow Z$ be the canonical  projection.
By (3) of prop. \ref{prop:automor-A0} and two lemmas above:  \notka{eq:deltaz-pi}
\begin{equation}\label{eq:deltaz-pi} (id\times \pi_Z^T\pi_Z)\delta_A\pi_Z^T=\delta_A\pi_Z^T.\end{equation}
Let the relation  $\delta_{Z}:Z\rel \Gamma_A\times Z$ be defined by \notka{def:tildez}
\begin{equation}\label{def:tildez} \delta_{Z}:=(id\times \pi_Z)\delta_A\pi_Z^T,\end{equation}
or in terms of graphs: $Gr(\delta_Z):=(id\times\pi_Z\times \pi_Z)Gr(\delta_A)$. This relation satisfies:
\notka{eq:coact-tildez}
\begin{equation}\label{eq:coact-tildez}
  (\delta_A\times id)\delta_{Z}=(id\times \delta_{Z})\delta_{Z}
  \end{equation}
Indeed:
\begin{equation*}
  \begin{split}
      (id\times \delta_{Z})\delta_{Z}&=(id\times (id\times \pi_Z)\delta_A\pi_Z^T)(id\times \pi_Z)\delta_A\pi_Z^T=
      (id\times id\times \pi_Z)(id\times \delta_A)(id\times \pi_Z^T\pi_Z)\delta_A\pi_Z^T=\\
      &=(id\times id\times \pi_Z)(id\times \delta_A)\delta_A\pi_Z^T=(id\times id\times \pi_Z)(\delta_A\times id)\delta_A\pi_Z^T=\\
      &=(\delta_A\times id) (id\times \pi_Z)\delta_A\pi_Z^T=(\delta_A\times id)\delta_{Z},
    \end{split}
  \end{equation*} 
  in the first equality definition of $\delta_{Z}$ i.e. (\ref{def:tildez}) was used, in the third one (\ref{eq:deltaz-pi}), in the next one 
  coassociativity of $\delta_A$  and finally (\ref{def:tildez}) again.\\
  In the next subsection, the canonical groupoid structure on $Z$ is described.

\subsection{Dividing a groupoid by a group action}

 The following  proposition describes division of a groupoid by an action of a group of automorphisms. \notka{prop:faktor}
 \begin{prop}\label{prop:faktor} Let $\Phi$ be a group of automorphisms of a groupoid $\Gamma$, i.e. each  $\varphi\in\Phi$ is  bijection
  $\Gamma\rightarrow \Gamma$ satisfying $\varphi(E)=E\,,\,\varphi s= s\varphi\,,\, \varphi m=m(\varphi\times\varphi)$.
  Let $\tilde{\Gamma}$ be the set of orbits of $\Phi$ and $\pi:\Gamma\rightarrow\tilde{\Gamma}$ the canonical projection.
  Assume that $\Phi$ satisfies:\notka{cond:free}
  \begin{equation}\label{cond:free}
    \forall\, \varphi\in\Phi, e\in E:\,\left[\varphi(e)=e\right]  \Rightarrow \left[ \forall\, \gamma\in F_l(e)\cup F_r(e): \varphi(\gamma)=\gamma\right]
    \end{equation}
    There exists unique structure $(\tilde{m},\tilde{s},\tilde{E})$ of a groupoid on $\tilde{\Gamma}$, such that $\pi^T:\tilde{\Gamma}\rel\Gamma$
    is a morphism.  Moreover $\pi^T$ is a monomorphism.
\end{prop} 
 \noindent Before the proof let us make some  comments.  First is that  the condition (\ref{cond:free}) 
excludes non trivial actions on groups, as it should be, since, in general, there is no group structure on the set of
orbits of a group of group automorphisms. The second is the lemma (the  proof is straightforward):
   
  \begin{lem} The following conditions are equivalent:\notka{lem:cond-free}
  \begin{enumerate}\label{lem:cond-free}
  \item Condition (\ref{cond:free}),
  \item $\displaystyle \forall\, \varphi\in\Phi, e\in E:\,\left[\varphi(e)=e\right]  \Rightarrow \left[ \forall\, \gamma\in F_l(e):
      \varphi(\gamma)=\gamma\right]$
  \item $\displaystyle \forall\, \varphi\in\Phi, e\in E:\,\left[\varphi(e)=e\right]  \Rightarrow \left[ \forall\, \gamma\in F_r(e):
      \varphi(\gamma)=\gamma\right]$
  \item  If $\gamma,\gamma'$ are in the same transitive component of $\Gamma$  then their isotropy groups (in $\Phi$) are equal.
  \end{enumerate} \komm{JESZCZE DLA PEWNOSCI SPRAWDZIC LEMAT--SPRAWDZONE!!!}
\end{lem}
\dowl

\noindent   The third one is that the action $R_a\,,\,a\in A_0$ in prop.\ref{prop:automor-A0} trivially fulfills the condition (\ref{cond:free}).
And the last  one is that some discussion is presented in Remarks after the proof.

\noindent {\em Proof:} In a   substantial part of the   proof  the condition (\ref{cond:free}) is not used.
It will be needed in  one (but crucial) step to show that one really obtains groupoid multiplication
and will be imposed at that moment. Let $h:\Gamma\rel\Gamma$ be the equivalence relation defined by the action of $\Phi$ i.e. orbit relation.
It is straightforward to check that this relation satisfies:
\begin{equation*}
     h(E)=E\,,\,hs=sh\,,\,hm\subset m(h\times h)
\end{equation*}
$\pi^T$ should be a morphism, therefore   relations 
$\tilde{m},\tilde{s}$ and the  subset $\tilde{E}$ must satisfy:\notka{piT-mor}
\begin{equation}\label{piT-mor}
  \pi^T(\tilde{E})=E\,,\,\pi^T \tilde{s}=s \pi^T\,,\,\pi^T \tilde{m}=m(\pi^T\times \pi^T).
  \end{equation}
  Because $\pi\pi^T=id_{\tilde{\Gamma}}$,  by applying $\pi$ to the equalities above one gets:
  \notka{def:tildestr}
\begin{equation}\label{def:tildestr} \tilde{E}=\pi(E)\,,\,\tilde{s}=\pi s \pi^T\,,\,\tilde{m}=\pi m(\pi^T\times \pi^T) \end{equation} 
so these relations are uniquely defined and it remains to check that they satisfy (\ref{piT-mor}) and define a groupoid structure.  
Notice that in terms of graphs, one has: \notka{tilde-graphG}  
\begin{equation}\label{tilde-graph}
Gr(\tilde{s})=(\pi\times\pi)Gr(s)\,,\, Gr(\tilde{m})=(\pi\times\pi\times\pi)Gr(m).
\end{equation}

First let us check (\ref{piT-mor}) is  satisfied. Let us compute:
$\displaystyle \pi^T(\tilde{E})=\pi^T\pi(E)=h(E)=E$; since $hs =s h$ one gets $\displaystyle \pi^T \tilde{s}=\pi^T \pi s \pi^T=s \pi^T$ by
lemma \ref{lem:hkomut}, and since $h m\subset m(h\times h)$, by the same lemma:   
$\displaystyle \pi^T \tilde{m}=\pi^T\pi m(\pi^T\times \pi^T)=m(\pi^T\times \pi^T)$.
Thus $(\tilde{m}, \tilde{s}, \tilde{E})$ defined by (\ref{def:tildestr}) satisfy
(\ref{piT-mor}).

For a groupoid structure on $\tilde{\Gamma}$ one needs (see \cite{SZ1}):
\begin{align}
  \tilde{m}(\tilde{m}\times id)& = \tilde{m}( id \times \tilde{m})\\
   \tilde{m}(\tilde{e}\times id)& = \tilde{m}( id \times \tilde{e})=id\\
  \tilde{s}^2 & =id\\
  \tilde{s}\tilde{m} & =\tilde{m}(\tilde{s}\times\tilde{s})\sim\\
  \forall \tilde{\gamma}\in \tilde{\Gamma}:\,\, & \emptyset \neq \tilde{m}(\tilde{s}(\tilde{\gamma}),\tilde{\gamma})\subset \tilde{E}\label{posit}
  \end{align}
  Let us compute:
\begin{equation*}\begin{split}
    \tilde{m}(\tilde{m}\times id)&=\pi m (\pi^T\times\pi^T)(\tilde{m}\times id)=
    \pi m (\pi^T \tilde{m} \times \pi^T)\stackrel{*}{=}\pi m (m(\pi^T\times\pi^T)\times \pi^T)=\\
    &=\pi m (m\times id) (\pi^T\times\pi^T\times \pi^T)\stackrel{**}{=}\pi m (id \times m) (\pi^T\times\pi^T\times \pi^T)=
    \pi m (\pi^T \times m(\pi^T\times \pi^T))\stackrel{***}{=}\\
    &=  \pi m (\pi^T \times \pi^T\tilde{m})=\pi m (\pi^T \times \pi^T)( id \times \tilde{m})=\tilde{m}(id\times\tilde{m}).
  \end{split}
\end{equation*} 
In (**) the associativity of $m$ was used and in (*) and (***) the equality: $\pi^T\tilde{m}=m (\pi^T\times \pi^T)$.
The next equality:
\begin{equation*}
    \tilde{m}(\tilde{e}\times id)=\pi m (\pi^T\times\pi^T)(\tilde{e}\times id)=\pi m (e\times \pi^T)=
    \pi m (e\times id) (id\times \pi^T)=\pi\pi^T=id,
\end{equation*}
and $\tilde{m}(id \times \tilde{e})$ can be computed in the same way.\\
Let us compute  $\tilde{s}^2$: $\displaystyle 
    \tilde{s}\tilde{s}=\pi s \pi^T \tilde{s}\stackrel{*}{=}\pi s s \pi^T=\pi \pi^T=id,$
    and in (*) the equality: $\pi^T \tilde{s}=s \pi^T$ was used. And the next one:
    \begin{equation*}\begin{split}
        \tilde{s} \tilde{m}&=\pi s \pi^T \tilde{m} =\pi s m(\pi^T\times \pi^T)=\pi m (s\times s)\sim(\pi^T\times \pi^T)=\\
        &=\pi m (s\pi^T \times s\pi^T )\sim=\pi m (\pi^T\tilde{s} \times\pi^T\tilde{s})\sim= \tilde{m}( \tilde{s} \times\tilde{s})\sim.
  \end{split}
\end{equation*}  
So up to now all conditions except possibly (\ref{posit}) are satisfied, and the assumption (\ref{cond:free}) was not used.

Multiplication relation $m$  in a groupoid  is a mapping from the subset $\Gamma^{(2)}\subset \Gamma\times \Gamma$ to $\Gamma$.
If one inspects the proof of this fact in \cite{SZ1}, one can see that the condition (\ref{posit}) is used. The present  proof shows that it is essential
and independent of other conditions. By taking $\Gamma=G$ - a group, and $\Phi$ the group of inner automorphisms of $G$ we obtain
$\tilde{\Gamma}$ as  the set of conjugacy classes in $G$ and the corresponding 
$\tilde{m}, \tilde{s}$ and $\tilde{e}$. Then it is easy to show an example that (\ref{posit}) is not satisfied (e.g. in $S_3$ group)
and $\tilde{m}$ is not a mapping.

It is clear (still without (\ref{posit})) that  for any $\tilde{\gamma}\in \tilde{\Gamma}$ the set
$\tilde{m}(\tilde{s}(\tilde{\gamma}),\tilde{\gamma})\neq \emptyset$, namely for any $\gamma$ with $\pi(\gamma)=\tilde{\gamma}$ one has
$(\pi(e_R(\gamma));\tilde{s}(\tilde{\gamma}), \tilde{\gamma}))\in \tilde{m}$. So we need:
$$(\tilde{\gamma}_1; \tilde{s}(\tilde{\gamma}), \tilde{\gamma})\in \tilde{m} \Rightarrow \tilde{\gamma}_1\in \tilde{E}$$
{\em Now, we assume (\ref{cond:free})} and compute:
$$(\tilde{\gamma}_1; \tilde{s}(\tilde{\gamma}), \tilde{\gamma})\in \tilde{m}\iff \exists \gamma_1, \gamma_2,\gamma\in\Gamma: \,
\pi(\gamma_1)=\tilde{\gamma}_1,\,,\,\pi(\gamma_2)=\tilde{s}( \tilde{\gamma})\,,\,\pi(\gamma)=\tilde{\gamma}\,{\rm and\,} (\gamma_1;\gamma_2,\gamma)\in m$$
Notice that $\pi(\gamma_2)=\tilde{s}( \tilde{\gamma})=\pi s \pi^T ( \tilde{\gamma})=\pi s \pi^T\pi (\gamma)= \pi s (\gamma)$;
this follows from $s^T=s$ and $\pi^T\pi s \pi^T=s \pi^T$ so also $\pi s=\pi s \pi^T \pi$. Thus $\gamma_2=\varphi(s (\gamma))$ for some $\varphi\in\Phi$ and
$(\gamma_1;\varphi(s (\gamma)),\gamma)\in m$. So $\el(\gamma)=\er(\varphi(s (\gamma)))=\varphi(\er(s (\gamma)))=\varphi(\el(\gamma))$ and by
(\ref{cond:free}): $\varphi(s (\gamma))=s (\gamma)$ therefore $\gamma_1=\er(\gamma)\in E$ and $\tilde{\gamma}_1\in \tilde{E}$ and (\ref{posit}) is proven.

The last thing is to prove that $\pi^T$ is a monomorphism. Monomorphisms in the Zakrzewski's category of groupoids were characterized in \cite{PS-short}.
It turns out that $h:\Gamma_1\rel\Gamma_2$ is a monomorphism iff {\em the kernel} of $h$ is $E_1$. The kernel is defined as
$ker(h):=\{\gamma_1\in\Gamma_1: (\gamma_2,\gamma_1)\in h\Rightarrow \gamma_2\in E_2\}$.
Now it is  clear that if $(e,\tilde{\gamma})\in \pi^T$ for some $e\in E$ then $\tilde{\gamma}=\pi(e)\in\tilde{E}$.
The proof is complete.
\dowl

\notka{rem:faktor}
\begin{re}\label{rem:faktor} There is the quotient construction of a groupoid $\Gamma\rra E$ by an action of a group of automorphisms $\Phi$
  \cite{Hig-Tayl,Tayl}. The resulting groupoid is denoted by $\Gamma\slash\hspace{-0.5ex}\slash \Phi$ and, in general,
  {\bf it is not a set of $\Phi$-orbits}.
  It is constructed as a quotient of a semi-direct product $\Gamma\rtimes  \Phi$ by the  normal subgroupoid $N$ generated by $E\rtimes \Phi$.
  Note, however, that unless normal subgroupoid $N\subset\Gamma$ is contained in the isotropy group bundle of $\Gamma$ the
  canonical map $\pi:\Gamma\rightarrow \Gamma\slash N$ is {\bf not} a (Zakrzewski) morphism. But condition (\ref{cond:free}) implies $E\rtimes \Phi$ {\bf is
  itself a normal subgroupoid} of $\Gamma\rtimes  \Phi$ and $(\Gamma\rtimes  \Phi)\slash (E\rtimes \Phi)$ can be identified with set
  of $\Phi$ orbits in $\Gamma$ and a groupoid structure is exactly as in Prop. \ref{prop:faktor}. Unless the action is trivial $E\rtimes \Phi$
  is not a group bundle and canonical projection is not a morphism.  
\end{re}
\notka{rem:faktor1}
\begin{re}\label{rem:faktor1} Prop. \ref{prop:faktor} can be thought of as a generalization of (algebraic part of) associated bundle construction.
  Let $H$ be a group and $X$ be a right, free $H$ space, and $\varphi$ a representation of $H$ on a vector space $V$.
  Consider the trivial transformation groupoid (group bundle) $\Gamma:X\rtimes V\rra X\times\{0\}$ with an action of
  $H$: $(x,v;h)\mapsto (xh, \varphi(h^{-1})(v))$. This is an action by automorphisms of $\Gamma$,
  due to the freeness of the action  condition (\ref{cond:free}) is satisfied
  and the resulting quotient is again a group bundle over $X\slash H$ (associated bundle).
  More general situation when the proposition can be applied is when
  $V$ is a group with a right action on $X$ and $\varphi: H\rightarrow Aut(V)$ is a homomorphism. Then the above action on
  $\Gamma$ is by automorphisms if $(xv)h=(xh)(\varphi(h^{-1})(v))$. Example of this more general situation is $X=V$ and $H\subset V$ is a subgroup.
  Both actions are by right translations and $\varphi(h)(v)=hvh^{-1}$.
\end{re}
\noindent Once  a groupoid structure on $Z$ is defined  one can prove that $\delta_Z$ is a morphism whenever $\delta_A$ is: 
\begin{prop} Let $(Z,m_Z,s_Z,E_Z)$ be $\Gamma_A\slash A_0$ with the structure defined in prop.\ref{prop:faktor} and
  $\delta_Z:Z\rel\Gamma_A\times Z$ be the relation  (\ref{def:tildez}). Then
  \begin{equation*} \delta_Z(E_Z)=A\times E_Z\,,\,\,\delta_Z s_Z=(s_A\times s_Z)\delta_Z\,,\,\,
    \delta_{Z} m_A\supset (m_A\times m_Z)(id\times\sim\times id)(\delta_{Z}\times \delta_{Z})
  \end{equation*}
  Moreover, if $\Gamma_A=G_A$ the last inclusion is  an equality.
  \end{prop}

  \noindent{\em Proof:} By (\ref{def:tildez}) $\delta_{Z}:=(id\times \pi_Z)\delta_A\pi_Z^T$ and by  prop.\ref{prop:faktor}
  $s_Z:=\pi_Z s_A \pi_Z^T,$ $\,m_Z=\pi_Z m_A(\pi_Z^T\times \pi_Z^T)$, and: 
\begin{equation*}
  \begin{split}(s_A\times s_Z)\delta_{Z}&=(s_A\times \pi_Z s_A \pi_Z^T)(id\times \pi_Z)\delta_A\pi_Z^T=
    (id\times\pi_Z) (s_A\times s_A )(id\times \pi_Z^T\pi_Z)\delta_A\pi_Z^T=\\
    &\stackrel{*}=(id\times\pi_Z) (s_A \times s_A )\delta_A\pi_Z^T
    \stackrel{**}=(id\times\pi_Z) \delta_A s_A\pi_Z^T\stackrel{***}=
    (id\times\pi_Z) \delta_A \pi_Z^Ts_Z=\delta_{Z}s_Z
  \end{split}
\end{equation*}
in $(*)$ the equality (\ref{eq:deltaz-pi}) was used, in $(**)$  the formulae (\ref{eq:amor-delta}), and $(***)$  is true since $\pi_Z^T$ is a morphism.

\noindent
For the multiplication:
\begin{equation*}
  \begin{split}
    \delta_{Z} m_Z&=(id\times \pi_Z)\delta_A\pi_Z^T m_Z\stackrel{*}=
    (id\times \pi_Z)\delta_A m_A(\pi_Z^T\times \pi_Z^T)\stackrel{**}\supset\\
    &\supset(id\times \pi_Z)(m_A\times m_A)(id\times\sim\times id) (\delta_A\times\delta_A) (\pi_Z^T\times \pi_Z^T)=\\
    &=(m_A\times \pi_Z m_A)(id\times\sim\times id) (\delta_A\pi_Z^T \times\delta_A\pi_Z^T)\stackrel{***}=\\
    &=
    (m_A\times \pi_Z m_A)(id\times\sim\times id)
    (id\times \pi_Z^T\times id \times \pi_Z^T) ((id\times \pi_Z)\delta_A\pi_Z^T \times(id\times \pi_Z)\delta_A\pi_Z^T)=\\
    &=(m_A\times \pi_Z m_A)(id\times id \times \pi_Z^T\times \pi_Z^T)(id\times\sim\times id)(\delta_{Z}\times \delta_{Z})=
  (m_A\times m_Z)(id\times\sim\times id)(\delta_{Z}\times \delta_{Z})
  \end{split}
\end{equation*}
$(*)$ holds since $\pi_Z^T$ is a morphism,  $(**)$ follows from (\ref{eq:amor-delta1}), $(***)$ follows from  (\ref{eq:deltaz-pi}),
and then definitions of $ \delta_{Z}$ and $m_Z$ were used. Notice also, that if $\Gamma_A=G_A$ (i.e. $(G;A,C)$ is DG) there is an equality in $(**)$
so the inclusion above is an equality as well.

\noindent Finally, for units:
  \begin{equation*}
  \begin{split}
    \delta_{Z}(E_Z) &=(id\times \pi_Z)\delta_A\pi_Z^T (E_Z)   =(id\times \pi_Z)\delta_A(E)=(id\times \pi_Z)(E\times E)=E\times E_Z,
  \end{split}
\end{equation*}
where the fact that $\pi_Z^T$ is a morphism and formulae (\ref{eq:amor-delta}) were used.

\dowl

To get the more explicit formula for $\delta_{Z}$ let us note that $c_L(\gamma)=c_L(\gamma a)$ for any $\gamma\in\Gamma$ and $a\in A_0$,
so there is a  well defined mapping $\tilde{c}_L:Z\ni [\gamma]\mapsto c_L(\gamma)\in C$; additionally if $c_R(\gamma)=\tilde{c}_L([\gamma_1])$ then
for any $\gamma'\in[\gamma_1]$: $\pi_Z(m_C(\gamma,\gamma_1))=\pi_Z(m_C(\gamma,\gamma'))\in Z$. By this remark and the definition of $\delta_A$
one gets $\delta_{Z}$:\notka{eq:tildeZ}
\begin{equation}\label{eq:tildeZ}
  \begin{split}\delta_{Z}:&=\{(\gamma, \pi_Z(m_C(s_C(\gamma),\gamma_1));z): \gamma,\gamma_1\in\Gamma\,,\, z\in Z\,,\,
    c_L(\gamma)=\tilde{c}_L(z)\,,\,\pi_Z(\gamma_1)=z\}=\\
    &=\{(\gamma, z;  \pi_Z(m_C(\gamma,\gamma_1))): \gamma,\gamma_1\in\Gamma\,,\,z\in Z\,,\,
    c_R(\gamma)=\tilde{c}_L(z)\,,\,\pi_Z(\gamma_1)=z\}
    \end{split}
  \end{equation}
At this moment we could  pass to the main example of this article, but before that let us briefly recall the classical situation.

\section{Oriented lines in Euclidean Space} 

Let $(M,V,\eta)$ be {\em an euclidean affine space} of dimension $\npj$,
where $\eta:V\times V\rightarrow \R$ is the euclidean scalar product;
the canonical isomorphism $V\rightarrow V^*$ will be denoted by $\eta$ as well.
Let $S^n:=\{p\in V: \eta(p,p)=1 \}\subset  V$.
A point  $m\in M$ and a vector $p\in S^n$ define the {\em oriented line} $l(p, m):=\{m+\lambda p: \lambda \in \R\}\subset M$.
The set of oriented lines in $M$ can be identified with the set of
classes of the equivalence relation on ``the sphere bundle'' of $M$ i.e. $S^n\times M$
\begin{equation}
  (p,m)\sim (p',m')\iff \left[ p=p'\,,\, m-m'=\lambda p\,,\,\lambda \in \R\right]
  \end{equation}
i.e. $[(p,m)]=l(p,m)$.

Let $ISO(\eta)$ be  the group of affine bijections of $M$ with linear part in $SO(\eta)$.
$ISO(\eta)$ acts transitively on $S^n\times M$ by $\tilde{\phi}(p,m):=(L_\phi(p), \phi(m))$,
where $L_\phi$ is the linear part of $\phi\in ISO(\eta)$;
it is straightforward to check that if $(p,m)\sim (p',m')$ then
$\tilde{\phi}(p,m)\sim \tilde{\phi}(p',m')$ so the action on $S^n\times M$ defines the action on lines (clearly transitive as well).

Since both actions are transitive, both spaces  can be identified with homogeneous spaces, after choosing some point
$(p_0,m_0)\in  S^n\times M$ or some line $l(p_0,m_0)$. The choice of $m_0$ enables us to identify $S^n\times M$ with 
$S^n\times V$ and $ISO(\eta)$ with $V\rtimes SO(\eta)$, and the choice of $p_0$ defines the  orthogonal decomposition
$V= p_0^\perp\oplus <p_0>=p_0^\perp\oplus\R$. With these identifications the isotropy group
of $(p_0, m_0)$ is $(0,SO(\eta^\perp))$ (where $\eta^\perp:=\eta\vert_{p_0^\perp}$ and
$SO(\eta^\perp)$ is included into $SO(\eta)$ as the stabilizer of $p_0$);
the isotropy group of the line $l(p_0,m_0)$ is $<p_0> \times SO(\eta^\perp)\subset V\rtimes SO(\eta)$
(one easily checks that this subgroup is a direct product of two of its subgroups as indicated by notation, since $p_0$ is fixed, it  is just
$\R \times SO(\eta^\perp)$). The set of lines in $M$ (or $V$ after $m_0$ was chosen)
can be identified with $TS^n=\{(p,v)\in S^n\times V: v\in p^\perp\}$ by $(p,v)\mapsto l(p, p+v)$. According to the decomposition 
$V=p_0^\perp \oplus <p_0>=p_0^\perp \oplus \R $  elements of $SO(\eta)$ can be written as
$\left(\begin{array}{cc} \Lambda, u \\ \eta(w), \alpha\end{array}\right)$,  where
$\Lambda:p_0^\perp\rightarrow p_0^\perp,\, u,w\in p_0^\perp,\, \alpha\in\R$.
The following maps identify $S^n\times V$ with $(V\rtimes SO(\eta))\slash (\{0\}\times SO(\eta^\perp))$ and
$TS^n$ with $(V\rtimes SO(\eta))\slash (<p_0> \times SO(\eta^\perp))$:
$$(V\rtimes SO(\eta))\slash (\{0\}\times SO(\eta^\perp))\ni [(v,A)]\mapsto
\left(p:=u+\alpha p_0 , v\right)\in S^n_\epsilon\times V$$
$$(V\rtimes SO(\eta))\slash (<p_0>\times SO(\eta^\perp))\ni [(v,A)]\mapsto
\left(p:=u+\alpha p_0, v-\eta(v,p)p \right)\in T S^n_\epsilon, $$
where $\displaystyle A:=\left(\begin{array}{cc} \Lambda, u \\ \eta(w), \alpha\end{array}\right)$. With these identifications, canonical left actions of
$V\rtimes SO(\eta)$ reads:\notka{eq:action-lines}
\begin{equation}\label{eq:action-lines}
  \begin{split}
(V\rtimes SO(\eta))\times (S^n\times V)\ni((z,A),(p,v))&\mapsto (Ap,z+Av)\in S^n\times V, \\
(V\rtimes SO(\eta))\times TS^n\ni((z,A),(p,v))& \mapsto (Ap,z+Av-\eta(z,Ap)Ap)\in TS^n
\end{split}
  \end{equation}
  Viewing $S^n\times V$ as  {\em the trivial bundle over $S^n$} one can identify it with
  $TS^n\oplus N$ (sum of tangent and (trivial) normal bundle of $S^n\subset V$) by
$$S^n\times V\ni (p,v)\mapsto (p; v-\eta(v,p)p, \eta(v,p)p)\in TS^n\oplus N$$
and the inverse mapping:
$$TS^n\oplus N\ni (p; v, \dot{s}p)\mapsto (p, v+\dot{s}p)\in S^n \times V$$
in this picture the action of $V\rtimes SO(\eta)$ reads:\notka{eq:action-spheres}
\begin{equation}\label{eq:action-spheres}
  ((z,A), (p; v, \dot{s}p))\mapsto (Ap;  z+Av-\eta(z,Ap)Ap,  (\dot{s}+\eta(z, Ap))Ap)
  \end{equation} 

  \section{The quantum space of euclidean lines}\label{sect:QE}

  The Quantum Euclidean Group, $QE(\npj)$ is defined by Double Lie Group coming from the Iwasawa decomposition of $SO_0(1,\npj)$.
  In the following I use the same notation as in \cite{PS-kappa2}.

\subsection{The Iwasawa  decomposition of $SO_0(1,\npj)$}\label{subs:Iwasawa}

\noindent
Let $G:=SO_0(1,\npj), n\geq 2$ and consider closed subgroups $B,C\subset G$:

{\em The group $B$} is $SO(\npj)$ embedded into $G$ by: $SO(\npj)\ni g\mapsto \left(\begin{array}{cc} 1 & 0 \\0 & g\end{array}\right)\in G$.
Elements of $SO(\npj)$ will be written as: \notka{def-B}
\begin{equation}\label{def-B}
(\Lambda,u,w,\alpha):=\left(\begin{array}{cc} 
\Lambda  & u\\
 w^t & \alpha
\end{array}\right),\quad\Lambda\in M_n(\R) , u,w\in\R^n\,,\,\alpha\in[-1,1],
\end{equation}
and $\Lambda, u, w, \alpha$ satisfy:  \notka{param-SOn}
\begin{equation}\label{param-SOn}\Lambda \Lambda^t+u u^t=I\,,\,\Lambda w+\alpha u=0\,,\,\Lambda^t \Lambda+ w w^t=I\,,\,
\Lambda^t u+\alpha w=0\,,\,|u|^2+\alpha^2=|w|^2+\alpha^2=1;
\end{equation}
these  equations imply that $\alpha=\det(\Lambda)$; $(\Lambda,u,w,\alpha)$ will denote also the corresponding element of $G$.

{\em The group $C$}  is: \notka{def-C}
\begin{equation}\label{def-C}
C:=\left\{\left(\begin{array}{ccc} 
\frac{s^2+1+|y|^2}{2 s} & -\frac1s y^t & \frac{s^2-1+|y|^2}{2 s}\\
-y & I & -y\\
\frac{s^2-1-|y|^2}{2 s} &  \frac1s y^t & \frac{s^2+1-|y|^2}{2 s}
\end{array}\right)\,s\in\R_+, y\in\R^n\right\}\subset G\,;
\end{equation}
it is isomorphic to the semi-direct product of $\R_+$ and $\R^n$  $\,\{(s,y)\in \R_+\times \R^n\}$  with  multiplication
$\displaystyle (s_1,y_1) (s_2,y_2):=(s_1 s_2, s_2 y_1+y_2)$.
As before  we will use $(s,y)$ to denote the corresponding element of $G$.

The Iwasawa decomposition for $G$ is $G=B C=C B$. To get its explicit form in the chosen parametrization one has to solve equation:
\notka{eq:iwasawa}
\begin{align}\label{eq:iwasawa} 
(\Lambda,u,w,\alpha)(s,y)&=(\tilde{s},\tilde{y})(\tilde{\Lambda}, \tilde{u}, \tilde{w}, \tilde{\alpha})
\end{align}
The solutions are presented in \cite{PS-kappa2}, Lemma 2.3 and the structure of the corresponding groupoid $G_B:G\rra B$ is described in
\cite{PS-kappa2}, Lemma 2.4.

It turns out that the set of one-point orbits of $G_B$ i.e. $N(C)\cap B=:B_0:=\{(\Lambda,0,0,1)\in B\}$ i.e. the group $SO(n)\subset SO(\npj)$
embedded by  $SO(n)\ni g\mapsto \left(\begin{array}{cc} g & 0 \\0 & 1\end{array}\right)\in SO(\npj)$, moreover: \notka{eq:iwasawa1}
\begin{equation}\label{eq:iwasawa1}
  (\Lambda,0,0,1)(s,y)=(s,\Lambda y )(\Lambda, 0, 0, 1)\,, \quad
  (\tilde{s},\tilde{y})(\tilde{\Lambda}, 0, 0, 1)=(\tilde{\Lambda}, 0, 0, 1) (\tilde{s},\tilde{\Lambda}^{-1} \tilde{y})
\end{equation}

\subsection{Quantum Euclidean Group}

The Quantum Group $QE(\npj)$ is obtained from DLG $(G;B,C)$ by the construction described in \cite{PS-DLG}. Its $C^*$-algebra is the (reduced)
$C^*$-algebra of a differential groupoid $G_B$ but, due to the special type of the group $C$, it can be said more.
A groupoid $G_B$ coming from the Iwasawa decomposition of $G$ is a transformation groupoid $G_B=B\rtimes C$ (the action is $(b,c)\mapsto b_R(bc)$).
Its reduced $C^*$-algebra $C^*_r(G_B)$ is the reduced crossed product $C(B)\rtimes_r C$ but since $C$ is amenable,
reduced and full crossed products coincide (\cite{DW}, Thm.7.13). Therefore for groupoid algebras we have $C^*(G_B)=C^*_r(G_B)$
(this is because, for transformation groupoids, the universal $C^*$-algebra of differential groupoid is ``between'' universal
and reduced crossed products, see \cite{PS-DG}).
Further, since $G_B\times G_B$ is also  the  transformation groupoid $G_B\times G_B=(B\times B)\rtimes (C\times C)$ by an amenable group $C\times C$,
again there is the equality $C^*(G_B\times G_B)=C^*_r(G_B\times G_B)$.
On the other hand $C^*_r(G_B\times G_B)=C^*_r(G_B)\otimes_{min}C^*_r(G_B)$ \cite{PS-DLG};
but  crossed product of nuclear $C^*$-algebras by amenable groups are nuclear (\cite{DW}, Cor. 7.18), therefore $C^*_r(G_B)$ is nuclear and there is only
one tensor product. The conclusion is: \notka{eq:tensor}
\newcommand{\tom}{\tilde{\om}}
\newcommand{\tof}{\tilde{f}}

\begin{equation}\label{eq:tensor}
  \begin{split}
  C^*(G_B\times G_B)=C^*_r(G_B\times G_B)=&C^*_r(G_B)\otimes_{min}C^*_r(G_B)=C^*_r(G_B)\otimes_{max}C^*_r(G_B)=\\=&C^*(G_B)\otimes_{min}C^*(G_B)=
  C^*(G_B)\otimes_{max}C^*(G_B).\end{split}
  \end{equation}

  As a consequence of the equality $C^*(G_B)=C^*_r(G_B)$ and equalities above the comultiplication for $QE(\npj)$ can be obtained as a direct $C^*$-lift
  of a differential groupoid morphism $\delta_B:=m_C^T: G_B\rel G_B\times G_B$ {\em without mediation of the multiplicative unitary};
  another direct application of this fact is that $SO(n)$ is a quantum subgroup of $QE(\npj)$ \cite{PS-future}.

  \subsection{The quantum space of euclidean lines} Now the construction described in section 2 will be performed for $G_B$ defined by the
  Iwasawa decomposition of $SO_0(1,\npj)$. The right action of $SO(n)=:B_0$ on $G_B$ is given by:
$$(\Lambda,u,w, \alpha; s, y)\Lambda_1=(\Lambda\Lambda_1,u,\Lambda_1^{-1}w, \alpha;s,\Lambda_1^{-1}y)\,,\quad \Lambda_1\in SO(n)$$

It is straightforward to verify that $Z:=G_B\slash B_0$ can be identified with $TS^n\times\R_+$ and the canonical projection $\pi_Z$ is given by
$$\pi_Z(\Lambda,u,w,\alpha; s, y):=
\left(p:=\left(\begin{array}{c} u\\ \alpha\end{array}\right),
  v:=\left(\begin{array}{c} \Lambda y\\ w^t y\end{array}\right)=:\left(\begin{array}{c} \mu\\ r\end{array}\right), s\right),$$
where $TS^n$ is identified with $\{(p,v)\in\R^{n+1}\times\R^{n+1}: |p|=1\,,\,p\cdot v=0\}$.
The next proposition describes the groupoid $Z\rra E_Z$. To simplify notation $(\tilde{m},\tilde{s}, \tilde{E})$ is used instead of $(m_Z, s_Z, E_Z)$.
\notka{prop:TS-deform}
\begin{prop}\label{prop:TS-deform} The groupoid $(Z,\tilde{m},\tilde{s},\tilde{E})$ has the following differential groupoid structure:
  \begin{itemize}
  \item[(a)] $Z=TS^n\times \R_+\,,\,\tilde{E}=S^n\times\{1\}\simeq S^n$.
  \item[(b)] Let $p_0:=\left(\begin{array}{c} {\bf 0} \\ 1\end{array}\right)\in S^n$, then left and right groupoid fibers over $p_0$ coincide and denoting them by
    $F_0$ we have $F_0=T_{p_0}S^n\times\R_+$; thus we have a decomposition $Z=F_0\sqcup Z_1$,
    where $F_0$  is a group and $Z$ is the closure of the open subgroupoid  $Z_1=T(S^n\setminus\{p_0\})\times \R_+$.
  \item[(c)] The group operations in  $F_0$ are:
    \begin{equation}
      (s_1, \mu_1) (s_2, \mu_2) =(s_1 s_2, s_2 \mu_1+\mu_2)\,,\,\, (s, \mu )^{-1}=\left(\frac1s, -\frac{\mu}{s}\right),
    \end{equation}
    where we write $(s, \mu)$ for $\left(p_0, \left(\begin{array}{c} \mu\\ 0\end{array}\right), s\right)\in F_0$; in other words $F_0\simeq C$
  \item[(d)] Let $x:=\frac{u}{1-\alpha}\in \R^n$ be the stereographic projection $S^n\setminus\{p_0\}\rightarrow \R^n$ and let $(x,\dot{x})$
    stands for the corresponding tangent coordinates: $T\R^n=\R^n\times \R^n\rightarrow T(S^n\setminus\{p_0\})$.
    The groupoid structure of $Z_1$ is given by:\notka{tildeZ-group}
    \begin{equation}\label{tildeZ-group}
      \begin{split}
        \tilde{s}(x,\dot{x},s)&=\left( s x -\frac{2}{1+|x|^2} \dot{x},\,
          -\frac{1+\left|s x - \frac{2\dot{x}}{1+|x|^2}\right|^2}{1+ |x|^2}\frac{\dot{x}}{s},\,\frac1s\right)\\
        \tilde{m}&=\left\{(x_1, s_2 \dot{x}_1+\dot{x}_2\frac{1+|x_1|^2}{1+|x_2|^2}, s_1 s_2;\,x_1, \dot{x}_1, s_1, x_2, \dot{x}_2, s_2):
          x_2=s_1 x_1- \frac{2 \dot{x}_1}{1+|x_1|^2}\right\}
      \end{split}
    \end{equation}
    these formulae imply that:  $\dsp \tilde{e}_L(x,\dot{x},s)=(x,0,1)\,,\,\,\tilde{e}_R(x,\dot{x},s)=(s x- \frac{2 \dot{x}}{1+|x|^2},0, 1)$.
  \item[(e)]  $Z_1$ is transitive and the isotropy group at $(x_0,0,1)$ is equal to $\{(x_0, x_0\frac{(s-1)(1+|x_0|^2)}{2}, s): s\in\R_+\}$.
  \item[(f)]  The closed submanifold  $\tilde{S}:=TS^n\times\{1\}\simeq TS^n$ is a wide subgroupoid of $Z$.
    
  \end{itemize}
\end{prop}
{\em Proof:\,\,} (a) The equality $\tilde{E}=S^n\times\{1\}$ follows from $\pi_Z(B)=S^n\times\{1\}$ and this is clear from the formula for $\pi_Z$,
and clearly this is a submanifold in $Z$.

Let us compute the groupoid inverse $\tilde{s}$. From prop.\ref{prop:faktor} $\tilde{s}$ is an involution $Z\rightarrow Z$
so it is a diffeomorphism if it is smooth.  The inverse in the groupoid $G_B$ is given by:
$$s_B(\Lambda,u,w,\alpha,s, y)=(\tilde{\Lambda}, \tilde{u}, \tilde{w},\tilde{\alpha}; \frac1s,-\frac ys),\quad{\rm where\,\,}
(\Lambda,u,w,\alpha,s, y)=(\tilde{s},\tilde{y})(\tilde{\Lambda}, \tilde{u}, \tilde{w},\tilde{\alpha}),$$
and, by (\ref{tilde-graph}),  for $\tilde{s}$ we have:
\begin{equation*}\begin{split}
    \tilde{s}(p,v,s)=(p_1, v_1, s_1)\iff \exists \, (\Lambda,u,w,\alpha,s, y):\, & \pi_Z(\Lambda,u,w,\alpha,s, y)=(p,v,s)\,,\,\\
    & \pi_Z(\tilde{\Lambda}, \tilde{u}, \tilde{w},\tilde{\alpha}; \frac1s,-\frac ys) =(p_1,v_1,s_1)
  \end{split}
\end{equation*}
Using solutions of equations (\ref{eq:iwasawa}) (see \cite{PS-kappa2}) one gets  $s_1=1\slash s$ and \notka{eq:inverse}
\begin{equation*}
  p_1=:\left( \begin{array}{c}u_1\\ \alpha_1\end{array} \right)\,,\,\, v_1=:\left( \begin{array}{c}\mu_1\\r_1\end{array}\right), {\rm \, where}\end{equation*}
\begin{equation}\label{eq:inverse}
  \begin{split}
u_1:= & \frac{2 (\frac{s-r}{1-\alpha} u -\mu)}{1+|\frac{s-r}{1-\alpha} u -\mu|^2}=
\frac{2 ((s-r)u  -(1-\alpha) \mu)}{(1-\alpha)(1+|\frac{s-r}{1-\alpha} u -\mu|^2)}=\frac{(s-r)u  -(1-\alpha) \mu}{Ms^2} \\
\alpha_1:=& \frac{|\frac{s-r}{1-\alpha} u -\mu|^2-1}{|\frac{s-r}{1-\alpha} u -\mu|^2+1}=
1-\frac{2(1-\alpha)}{(1-\alpha)(1+ |\frac{s-r}{1-\alpha} u -\mu|^2)}=1-\frac{1-\alpha}{Ms^2}\\
\mu_1:=& -\frac{r_1(s-r)s +r}{(1-\alpha)s } u +(r_1-\frac1s)\mu= -\frac{(2s-r)(|\mu|^2+r^2)+r(1-s^2)}{2Ms^3} u +(r_1-\frac1s)\mu\\
r_1:=&- \frac{1}{ M s^3 (1-\alpha)}\left[(s-r) u -(1-\alpha) \mu\right]^t\left[ r u +(1-\alpha)\mu\right]=\frac{(1-\alpha)(|\mu|^2+r^2)-rs}{M s^3}
\end{split}
\end{equation}
and 
\begin{equation*}
  M:=\frac{(1-\alpha)}{2 s^2}\left(1+ |\frac{s-r}{1-\alpha} u -\mu|^2\right)=\frac{1}{2 s^2}\left[ (1-\alpha)(|\mu|^2+r^2+1)+ (1+\alpha) s^2 -2rs\right]\,,\,\,
  \end{equation*}
  From the first expression for $M$ it is clear that $M> 0$ for $\alpha\neq 1$ and, from the second one, $M=1$ for $p=p_0$ i.e $\alpha=1, u=0$
  (then $r=0$) so $M$ is strictly positive everywhere, and formulae (\ref{eq:inverse}) are well defined on $Z$;
  short examination proves they define smooth mapping of $Z$. For $p=p_0$ we get:
 \begin{equation}\label{tildes-p0} \tilde{s}\left( p_0, v, s\right)=
  \left( p_0, -\frac{v}{s}, \frac{1}{s}\right),
\end{equation}
so $T_{p_0}S^n\times\R_+$ and $Z_1$ are $\tilde{s}$ invariant.
  
The stereographic coordinates on $S^n\setminus\{p_0\}$ are given by:
$$\R^n\ni x\mapsto p(x)=\left(\begin{array}{c}u=\frac{2 x}{1+|x|^2}\\\alpha=\frac{|x|^2-1}{|x|^2 +1}\end{array}\right)\in
S^n\setminus\{p_0\}.$$
Let $(x,\dot{x})$ be coordinates on $\R^n\times \R^n=T\R^n$: $(x,\dot{x})$ represents the vector tangent at $t=0$ to the
curve $t\mapsto x+t \dot{x}$;  we get coordinates on $T(S^n\setminus\{p_0\})$:
\notka{eq:stereo}
\begin{equation}\label{eq:stereo}
  \begin{split}
    (x,\dot{x})& \mapsto \left (p(x) , v(x,\dot{x})=\left(\begin{array}{c}\mu(x,\dot{x}) \\ r(x,\dot{x})\end{array}\right)\right) \\
    \mu(x,\dot{x})&:=-\frac{4 x^t\dot{x}}{(1+|x|^2)^2} x+\frac{2 \dot{x}}{1+|x|^2}\,,\,r(x,\dot{x}):=\frac{4 x^t\dot{x}}{(1+|x|^2)^2},
  \end{split}
  \end{equation}
and the inverse mapping:\notka{eq:stereo-inv}
\begin{equation}\label{eq:stereo-inv}
  \left( p=\left(\begin{array}{c}u\\ \alpha\end{array}\right),v=\left(\begin{array}{c}\mu \\ r\end{array}\right)\right) \mapsto
  \left( x:=\frac{u}{1-\alpha},\, \dot{x}:=\frac{r}{(1-\alpha)^2}u+ \frac{1}{1-\alpha}\mu\right)
  \end{equation}
Using these formulae and (\ref{eq:inverse}) one gets formulae for $\tilde{s}$ on $Z_1$:\notka{tildes-Z1}
\begin{equation}\label{tildes-Z1}
  \tilde{s}(x,\dot{x},s)=\left( s x -\frac{2}{1+|x|^2} \dot{x},\,
    -\frac{1+\left|s x - \frac{2\dot{x}}{1+|x|^2}\right|^2}{1+ |x|^2}\frac{\dot{x}}{s},\,\frac1s\right)
\end{equation}

Let us now find groupoid projections; for a groupoid $\Gamma\rra E$ we have in general $e=e_L(\gamma)\,\iff\, (\gamma;e, \gamma)\in m$.
Therefore
\begin{equation*}
  \begin{split} &(p_1,0,1) =\tilde{e}_L(p,v,s) \iff \left(p,v,s;p_1,0,1;p,v,s\right)\in \tilde{m}\iff \\
    & \exists\, b_1,b_2,b_3\in B, c_1,c_3\in C:
    \pi_Z(b_1 c_1)=(p,v,s)=\pi_Z(b_3 c_3)\,,\,\pi_Z(b_2)=(p_1,0,1)\,,\, (b_1 c_1;b_2, b_3c_3)\in m_B,
  \end{split}
\end{equation*}
and it follows  $b_1=b_2=b_3$,  $c_1=c_3$, and $p_1=p$. This way we get $\tilde{e}_L(p,v,s)=(p,0,1)\in \tilde{E}$, so this is clearly a submersion and since
$\tilde{s}$ is a diffeomorphism and $\tilde{e}_R=\tilde{e}_L\cdot \tilde{s}$, $\tilde{e}_R$ is a submersion as well. 
In particular for  $p=p_0$ using (\ref{tildes-p0}):
$$\tilde{e}_R(p_0,v,s)=\tilde{e}_L(\tilde{s}(p_0,v,s))=(p_0,0,1)=\tilde{e}_L(p_0,v,s),$$
what proves (b); the formula for inverse in $F_0$ follows directly from  (\ref{tildes-p0}).
To get a formula for multiplication in $F_0$ let us write elements from $F_0$ as $(s,\mu)$ and observe that
$(s_3,\mu_3)=\tilde{m}(s_2,\mu_2, s_1,\mu_1)$ iff exist $b_1, b_2, b_3\in B_0$ and $c_1, c_2, c_3\in C$ with
$\pi_Z(b_i c_i)=(s_i, \mu_i)$ and $b_3 c_3=m_B(b_2 c_2, b_1 c_1)$. Since $b_2\in B_0$, it follows $b_1=b_2=b_3=:\Lambda$ and $c_3=c_2 c_1$ so
$(s_3, y_3)=(s_1 s_2, s_2 y_1+y_2)$, and $\mu_3=\Lambda y_3=\Lambda(s_2 y_1+y_2)=s_2 \mu_1+\mu_2$. This way (c) is proven.

The last thing is to prove the formula for multiplication in $Z_1$.
Notice that if $b_3 c_3=m_B(b_2 c_2, b_1 c_1)$ then $c_3=c_2 c_1$ and by the formula for $\pi_Z$ we obtain
$$(p_3, v_3, s_3)=\tilde{m}(p_1, v_1, s_1; p_2, v_2, s_2)\Rightarrow s_3=s_1 s_2.$$
By the formula for $\tilde{e}_L$ and (\ref{tildes-Z1}) one gets
$\tilde{e}_R$ on $Z_1$:
$$\tilde{e}_R(x,\dot{x},s)= (s x -\frac{2}{1+|x|^2} \dot{x},0, 1).$$
For multiplication relation observe that necessary and sufficient condition for $(x_1,\dot{x}_1,s_1)$ and
$(x_2,\dot{x}_2,s_2)$ being composable is $\tilde{e}_R(x_1,\dot{x}_1,s_1)=\tilde{e}_L(x_2,\dot{x}_2,s_2)$ i.e.
$x_2=s_1 x_1 -\frac{2}{1+|x_1|^2} \dot{x}_1$. Assuming this is the case and comparing left and write projections of  product we obtain:
$$(x_3,\dot{x}_3,s_3)=\tilde{m}( x_1,\dot{x}_1,s_1, x_2,\dot{x}_2,s_2)\Rightarrow
x_3=x_1\,,\, \tilde{e}_R(x_3,\dot{x}_3,s_3)=\tilde{e}_R(x_2,\dot{x}_2,s_2),$$
and
$$s_3 x_1 -\frac{2}{1+|x_1|^2} \dot{x}_3=s_2(s_1 x_1 -\frac{2}{1+|x_1|^2} \dot{x}_1)-\frac{2}{1+|x_2|^2} \dot{x}_2 .$$
Since $s_3=s_2 s_1$ we obtain:
$$\dot{x}_3=s_2  \dot{x}_1+ \frac{1+|x_1|^2}{1+|x_2|^2} \dot{x}_2.$$
This way (d) is proven.  Since $\tilde{e}_R$ and $\tilde{e}_L$ are submersions, the set of composable elements in $Z$ is a submanifold;
from formula for $\tilde{m}$ (\ref{tildeZ-group})  it is clear that it is smooth on (composable elements from ) $Z_1\times Z_1$, so it remains to
prove that $\tilde{m}$ is smooth at points of $F_0$. 

Let  $S^n_+:=\left\{\left(\begin{array}{c}u\\ \alpha\end{array}\right)\in S^n: \alpha>0\right\}$, and let us use $(u,\dot{u},s)$ as coordinates on $TS^n_+\times \R_+$ i.e.
\begin{align*}
U&:=\{(u,\dot{u},s)\in \R^n\times\R^n\times \R_+ :|u|< 1\}\\
\Phi: & \, U \ni (u,\dot{u},s)\mapsto
\left( p=\left(\begin{array}{c}u\\ \alpha(u)\end{array}\right),
  v=\left(\begin{array}{c}\dot{u} \\ \frac{- u^t \dot{u}}{\alpha(u)}\end{array}\right), s \right)\in TS^n_+\times \R_+\,,
 \,\, \alpha(u):=\sqrt{1-|u|^2}.
\end{align*}
To find the formula for $\tilde{m}$ in coordinates $(u,\dot{u},s)$ we begin by finding $\tilde{e}_R$. The set
$W:=\tilde{e}_R^{-1}(S^n_+)\cap \tilde{e}_L^{-1}(S^n_+)$ is an open neighborhood of $F_0$. Using the formula for $\tilde{e}_R$ in stereographic coordinates
one gets for $u\neq 0$ that
$$(u,\dot{u},s)\in \Phi^{-1}(W)\,\,\iff\,\, \left|(s+\frac{1}{\alpha}u^t \dot{u}) u -(1-\alpha)\dot{u}\right|> (1-\alpha)\,,\,\alpha:=\alpha(u).$$
By squaring the inequality above we obtain the equivalent (for $u\neq 0$)  inequality:\notka{def:A}
\begin{equation}\label{def:A}
  A(u,\dot{u},s):=(s+\frac{1}{\alpha}u^t \dot{u})^2 (1+\alpha)-2 (s+\frac{1}{\alpha}u^t \dot{u})u^t \dot{u}+(1-\alpha)(|\dot{u}|^2-1)>0,
  \end{equation}
  which is satisfied also for $u=0$; this way $\Phi^{-1}(W)=\{(u,\dot{u},s)\in U: A(u,\dot{u},s)>0\}$, and  for $(u,\dot{u},s)\in \Phi^{-1}(W)$:
  \notka{eq:er-u}
  \begin{equation}\label{eq:er-u}
    \tilde{e}_R(u,\dot{u},s)=(\tilde{u},0,1)\,,\,
    \tilde{u}:=\frac{2}{A(u,\dot{u},s)+2(1-\alpha)}\left((s+\frac{1}{\alpha}u^t \dot{u}) u -(1-\alpha)\dot{u}\right)
  \end{equation}
  To simplify notation we will write:
  \begin{equation*}
    \tilde{e}_R(u,\dot{u},s) =\tilde{u}(u,\dot{u},s)\,\,,\,\,k(u,\dot{u},s):=(s+\frac{1}{\alpha}u^t \dot{u}) u -(1-\alpha)\dot{u}\,,\quad{\rm\,then}\quad
                          \tilde{u}:=\frac{2k(u,\dot{u},s)}{A(u,\dot{u},s)+2(1-\alpha)}.
  \end{equation*}
  Notice that 
  $\dsp \left|(s+\frac{1}{\alpha}u^t \dot{u}) u -(1-\alpha)\dot{u}\right|^2=|k(u,\dot{u},s)|^2=(1-\alpha)(A(u,\dot{u},s)+(1-\alpha)).$
  
  By straightforward computation one gets that for fixed $(u,s)$ with $0<|u|<1$ the mapping:
  $$\{\dot{u}: (u,\dot{u}, s)\in \Phi^{-1}(W)\}\ni \dot{u}\mapsto  \tilde{e}_R(u,\dot{u},s)
=:\tilde{u}=\frac{2k(u,\dot{u},s)}{A(u,\dot{u},s)+2(1-\alpha)}\in \{\tilde{u}\in\R^n: 0<|\tilde{u}|<1\}$$
is invertible. Its inverse is given by:\notka{eq:er-u-inv}
\begin{equation}\label{eq:er-u-inv}
\{\tilde{u}\in\R^n: 0<|\tilde{u}|<1\}\ni \tilde{u}\mapsto \dot{u}:=
\left(\frac{u^t \tilde{u}}{1-\tilde{\alpha}}-\alpha s\right)\frac{u}{1-\alpha}-\frac{\tilde{u}}{1-\tilde{\alpha}}\in
\{\dot{u}: (u,\dot{u}, s)\in \Phi^{-1}(W)\}
\end{equation}
and the following equalities hold:\notka{eq:alpha-A}
\begin{equation}\label{eq:alpha-A}
  1-\alpha(\tilde{u})=\frac{2(1-\alpha(u))}{A(u,\dot{u},s)+2(1-\alpha(u))}\,\,,\,\,\,
  \frac{\tilde{u}}{1-\alpha(\tilde{u})}=\frac{k(u,\dot{u},s)}{1-\alpha(u)}
    \end{equation}
  Let $z_1^0:=(p_0, v_1^0, s_1^0)=\Phi(0,\dot{u}_1^0, s_1^0)$ and $z_2^0:=(p_0, v_2^0, s_2^0)=\Phi(0,\dot{u}_2^0, s_2^0)$. Let $\sO\subset \R^n$ be an
  open  neighborhood of $0$ and $\sO_2\subset \R^n\times\R_+$ be an open  neighborhood of $(\dot{u}_2^0, s_2^0)$ such that
  $\Phi(\sO\times\sO_2)\subset W$; moreover let $\sO' \subset \R^n$ be an
  open  neighborhood of $0$ and $\sO_1\subset \R^n\times\R_+$ be an open  neighborhood of $(\dot{u}_1^0, s_1^0)$ such that
  $\tilde{e}_R(\Phi(\sO'\times\sO_1))\subset \sO$, then
  $$\sO'\times\sO_1\times \sO_2\ni(u_1,\dot{u}_1, s_1,\dot{u}_2,s_2)\mapsto (\Phi(u_1,\dot{u}_1, s_1), \Phi(\tilde{u}_1,\dot{u}_2, s_2))\in Z^{(2)}\,,\,
  (\tilde{u}_1,0,1):=\tilde{e}_R(u_1,\dot{u}_1,s_1)$$
  are local coordinates in $Z^{(2)}$ around $(z_1^0, z_2^0)$. In this coordinates:
  $$\tilde{m}(u_1,\dot{u}_1, s_1,\dot{u}_2,s_2)=(u_1, \dot{u}_3, s_1 s_2)$$
  For $u_1=0$ we already have that  $\dot{u}_3=s_2\dot{u}_1+\dot{u}_2$, so it remains to compute $\dot{u}_3$ for $u_1\neq 0$ i.e. translate
  (\ref{tildeZ-group}) into $(u,\dot{u})$ coordinates. Let
  $\tilde{u}_1:=\tilde{e}_R(u_1,\dot{u}_1, s_1)=\frac{2 k_1}{A_1+2(1-\alpha_1)}$ and
  $\tilde{u}_2:=\tilde{e}_R(\tilde{u}_1,\dot{u}_2,s_2)=\frac{2 k_2}{A_2+2(1-\tilde{\alpha}_1)}$, where $k_1:=k(u_1,\dot{u}_1,s_1)$,
  $A_1:=A(u_1,\dot{u}_1,s_1)$, and similarly for $k_2$ and $A_2$. Using  the equality
  $\tilde{u}_2=\tilde{e}_R(u_1, \dot{u}_3, s_1 s_2)$, denoting $\tilde{\alpha}_2:=\alpha(\tilde{u}_2)$, by  (\ref{eq:er-u-inv}) and
  (\ref{eq:alpha-A}), one gets:
  \begin{equation*}\begin{split}
  \dot{u}_3 & =\left(\frac{u_1^t\tilde{u}_2}{1-\tilde{\alpha}_2}-\alpha_1 s_1 s_2\right) \frac{u_1}{1-\alpha_1}-\frac{\tilde{u}_2}{1-\tilde{\alpha}_2}=
  \left(\frac{u_1^t k_2}{1-\tilde{\alpha}_1}-\alpha_1 s_1 s_2\right) \frac{u_1}{1-\alpha_1}-\frac{k_2}{1-\tilde{\alpha}_1}=\\
  &= \left[\left(s_2+\frac{\tilde{u}_1^t\dot{u}_2}{\tilde{\alpha}_1}\right)\frac{u_1^t\tilde{u}_1}{1-\tilde{\alpha}_1}-u_1^t\dot{u}_2-\alpha_1 s_1 s_2\right]
  \frac{u_1}{1-\alpha_1}- \left(s_2+\frac{\tilde{u}_1^t\dot{u}_2}{\tilde{\alpha}_1}\right)\frac{\tilde{u}_1}{1-\tilde{\alpha}_1}+\dot{u}_2\end{split}
\end{equation*}
In the last equality the definition of $k_2$ was used:
$\dsp k_2=(s_2+\frac{\tilde{u}_1^t\dot{u}_2}{\tilde{\alpha}_1})\tilde{u}_1-(1-\tilde{\alpha}_1) \dot{u}_2$. By (\ref{eq:alpha-A}) again 
$\frac{\tilde{u}_1}{1-\tilde{\alpha}_1}=\frac{k_1}{1-\alpha_1}$, and substituting $k_1=(s_1+\frac{u_1^t\dot{u}_1}{\alpha_1})\frac{u_1}{1-\alpha_1}-(1-\alpha_1) \dot{u}_1$:
$$\dot{u}_3=K \frac{u_1}{1-\alpha_1}+ (s_2+\frac{\tilde{u}_1^t\dot{u}_2}{\tilde{\alpha}_1})\dot{u}_1+\dot{u}_2,$$
where
$$K:=(s_2+\frac{\tilde{u}_1^t\dot{u}_2}{\tilde{\alpha}_1})\frac{u_1^t k_1}{1-\alpha_1}-u_1^t\dot{u}_2-\alpha_1 s_1 s_2 -
(s_2+\frac{\tilde{u}_1^t\dot{u}_2}{\tilde{\alpha}_1})(s_1+\frac{u_1^t\dot{u}_1}{\alpha_1}),$$
and, after short computation:
\begin{equation*}
  \begin{split}
    K&=-(1-\alpha_1)\left[ \left(s_1+\frac{u_1^t\dot{u}_1}{\alpha_1}\right)^2+ |\dot{u}_1|^2-1\right]\\
    s_2+\frac{\tilde{u}_1^t\dot{u}_2}{\tilde{\alpha}_1}&=s_2+\frac{2 u_1^t\dot{u}_2}{A_1}\left(s_1+\frac{u_1^t\dot{u}_1}{\alpha_1}-(1-\alpha_1)\right),
  \end{split}\end{equation*}  and finally:\notka{eq:dotu3}
\begin{equation}\label{eq:dotu3}
  \dot{u}_3=- \left[ \left(s_1+\frac{u_1^t\dot{u}_1}{\alpha_1}\right)^2+ |\dot{u}_1|^2-1\right] u_1+
  \left[s_2+\frac{2 u_1^t\dot{u}_2}{A_1}\left(s_1+\frac{u_1^t\dot{u}_1}{\alpha_1}-(1-\alpha_1)\right)\right] \dot{u}_1+\dot{u}_2
\end{equation}
This expression is well defined also for $u_1=0$ (and it gives $\dot{u}_3=s_2\dot{u}_1+\dot{u}_2$ as it should) and is
smooth on $\sO'\times\sO_1\times \sO_2$.  This completes the proof that $Z$ is a differential groupoid.

e) and (f) are  direct consequences of (a)-(d).
\dowl

Now, structure of the groupoid $\tilde{S}:TS^n\rra S^n$ will be described explicitly.
$\tilde{S}$ is a disjoint union of $T_{p_0}S^n$ (this is the isotropy group of $p_0$) with the standard (vector)  group structure and open,
transitive subgroupoid $\tilde{S}_1:=T(S^n\setminus\{p_0\})$ with no isotropy (in other words $\tilde{S}_1$
is isomorphic to the pair groupoid $\R^n\times\R^n$ by the map $(x, \dot{x})\mapsto (x, e_R(x,\dot{x})$)
and the following groupoid structure in stereographic coordinates:
\notka{tildeT-group}
\begin{equation}\label{tildeT-group}
      \begin{split}
        \tilde{s}(x,\dot{x})&=\left( x -\frac{2}{1+|x|^2} \dot{x},\,
          -\frac{1+\left| x - \frac{2\dot{x}}{1+|x|^2}\right|^2}{1+ |x|^2}\dot{x}\right)\\
        \tilde{m}&=\left\{(x_1, \dot{x}_1+\dot{x}_2\frac{1+|x_1|^2}{1+|x_2|^2};\,x_1, \dot{x}_1, x_2, \dot{x}_2):
          x_2=x_1- \frac{2 \dot{x}_1}{1+|x_1|^2}\right\}
      \end{split}
    \end{equation}

    The left projection in $\tilde{S}$ is the bundle projection, and the right one, in stereographic coordinates, is given by
    $e_R(x,\dot{x})=(x- \frac{2 \dot{x}}{1+|x|^2}, 0)$. To get expression for right projection and multiplication in  $(u,\dot{u})$ coordinates
    it is enough to put $s_1=s_2=1$ in (\ref{eq:dotu3}) and (\ref{eq:er-u}).
    
\subsection{The Lie algebroid of $Z$. }  The left projection in $Z=TS^n\times \R_+$ is $(p,v,s)\mapsto (p,0,1)$, 
therefore the Lie algebroid of $Z$ (the bundle of vectors tangent to left fibers at $\tilde{E}$), denoted by $\lie{Z}$,
is just the bundle $TS^n\oplus \bunR$, where $\bunR$ is the trivial bundle $S^n\times \R$.
An  element $(p,\dot{v},\dot{s})\in  TS^n\oplus \bunR$ as  an element of  $T_{(p,0,1)}(TS^n\times \R_+)$
is represented by  the  curve $]-\epsilon,\epsilon[\ni t\mapsto (p, \dot{v}t, 1+\dot{s}t)\in Z$.  
\komm{ZMIANA KONWENCJI: TO CO BYLO $X_0$ JEST TERAZ $X_{n+1}$}
As a vector bundle $\lie{Z}=TS^n\oplus \bunR$ is trivial; it can be identified with $S^n\times \R^{n+1}$ by 
$TS^n\oplus \bunR\ni (p,\dot{v},\dot{s})\mapsto (p,\dot{v}+\dot{s}p)\in S^n\times \R^{n+1}$
(i.e. $\bunR=S^n\times \R$ is identified with the normal bundle to $S^n$ in $\R^{n+1}$);
the inverse mapping is: $S^n\times \R^{n+1}\ni (p,w)\mapsto (p, w-(\skalp{w}{p})p,  \skalp{w}{p})\in TS^n\oplus \bunR$.
For $p\in S^n$ let the  projection $\R^{n+1}\ni w\mapsto w-(\skalp{w}{p}) p \in p^\perp\subset \R^{n+1}$ (considered as $T_pS^n$) be denoted by $P_p(w)$.
Taking the canonical basis $e_i$ in $\R^{\npji}$ we obtain the
basis of sections of $\lie{Z}$:\notka{def:basis-sections}
\begin{equation} \label{def:basis-sections} \tilde{X}_i(p,0,1):=(p, \dot{v}_i(p), \dot{s}_i(p))\,,\,\,\,\,
                    \dot{v}_i(p):=P_p(e_i)=e_i-p^i p\,,\, \dot{s}_i(p):=p^i, \,\,i=1,\dots,\npj
\end{equation}
In particular,  $\dot{v}_i(p_0)=(1-\delta_{i, n+1}) e_i\,,\,  \dot{s}_i(p_0)=\delta_{i, n+1}$.\notka{prop:liealg}
\begin{prop} \label{prop:liealg} The Lie algebroid structure of $\lie{Z}$ is given by the bracket:\notka{eq:liealg-nawias}
  \begin{equation}\label{eq:liealg-nawias}\begin{split}
      [ \tilde{X}_i,  \tilde{X}_j]&=2 (p^j \tilde{X}_i-p^i \tilde{X}_j)\,,\,\,\,i,j=1,\dots,  n\\
      [\tilde{X}_i,  \tilde{X}_{n+1}]&=2(p^{n+1} \tilde{X}_i- p^i \tilde{X}_{n+1}) -\tilde{X}_i+2 p^i\sum_{k=1}^{n+1}p^k \tilde{X}_k\,,\,\,\,i=1,\dots,  n
    \end{split}\end{equation} 
  and the anchor:\notka{eq:liealg-kotwa}
  \begin{equation} \label{eq:liealg-kotwa}
    \rho(\tilde{X}_i)(p)=-(1-p^{n+1})(e_i-p^ip)+ p^i(e_{n+1}-p^{n+1}p)=P_p(p^ie_{n+1}+p^{n+1}e_i-e_i)\,,\,\,i=1,\dots,  n+1 \end{equation}     
\end{prop} 
\begin{re}
  By identification with the trivial bundle $S^n\times \R^{n+1}$, this can be interpreted as an algebroid structure on this bundle. Using
  constant sections $\tilde{e}_i(p):=e_i$, the bracket and the anchor can be rewritten as:
  \begin{equation}
    \begin{split}
      [\tilde{e}_i, \tilde{e}_j](p)=&=2 (p^j \tilde{e}_i(p) -p^i \tilde{e}_j(p))\,,\,\,\,i,j=1,\dots,  n\\
      [\tilde{e}_i,  \tilde{e}_{n+1}](p) &=2(p^{n+1} \tilde{e}_i(p) - p^i \tilde{e}_{n+1}(p) ) +(I-2P_p)(e_i)\\
      \rho(p,w)&:=-(1-p^{n+1})P_p(w)+(\skalp{p}{w})P_p(e_{n+1})=P_p((\skalp{p}{e_{n+1}})w+(\skalp{p}{w})e_{n+1}-w) 
    \end{split}
  \end{equation}
\end{re}

  \noindent {\em Proof: } First we compute brackets in stereographic coordinates (\ref{eq:stereo},\ref{eq:stereo-inv}).  The tangent space to the left fiber
  at $(x,0,1)$ is spanned by linearly independent vectors: 
  \begin{equation}
    X^l_{n+1}(x,0,1):=\partial_s  \,\,{\rm and}\,\,\,  X^l_i(x,0,1):=\partial_{\dot{x}^i}\,,\,\, i=1,\dots,n.
  \end{equation} 
    Using the formula for multiplication (\ref{tildeZ-group}) one obtains corresponding left invariant vector fields on $Z_1$ 
    ( $\tilde{x}:= s x -\frac{2}{1+|x|^2} \dot{x}$ i.e. $\tilde{e}_R(x,\dot{x},s)=(\tilde{x},0,1)$, cf. (\ref{tildeZ-group}) ): 
\begin{equation} 
  X_{n+1}^l(x,\dot{x},s):=\sum_{k=1}^{n}\dot{x}^k\partial_{\dot{x}^k}+ s\partial_s\,,\,\,
  X_i^l(x,\dot{x},s):=\frac{1+|x|^2}{1+|\tilde{x}|^2} \partial_{\dot{x}^i}
\end{equation} 
Having these fields, by standard computations, one gets: 
\begin{equation}
  [X_i^l, X_j^l]=\frac{4}{(1+|\tilde{x}|^2)^2}\left(\tilde{x}^i X_j^l-\tilde{x}^j X_i^l\right)\,,\,
  [X_i^l, X_{n+1}^l]=(3-\frac{2}{1+|\tilde{x}|^2}) X_i^l\,,\,i,j=1,\dots,n
  \end{equation} 
The extension of sections (\ref{def:basis-sections}) to left invariant vector fields (in stereographic coordinates)  reads: 
\begin{equation}\begin{split}
      \tilde{X}_i(x,\dot{x},s)&:=\frac{2  \tilde{x}^i}{1+|\tilde{x}|^2}\left(\sum_{k=1}^n \dot{x}^k\partial_{\dot{x}^k}+s\partial_s\right)
      -\frac{1+|x|^2}{1+|\tilde{x}|^2} \,\tilde{x}^i \sum_{k=1}^n \tilde{x}^k\partial_{\dot{x}^k} +\frac{1+|x|^2}{2}\partial_{\dot{x}^i}\,,\,i=1,\dots,n\\
    \tilde{X}_{n+1}(x,\dot{x},s)&:=\frac{|\tilde{x}|^2-1}{1+|\tilde{x}|^2} \left(\sum_{k=1}^n \dot{x}^k\partial_{\dot{x}^k}+s\partial_s\right)+
    \frac{1+|x|^2}{1+|\tilde{x}|^2} \, \sum_{k=1}^n \tilde{x}^k\partial_{\dot{x}^k}
    \end{split}
\end{equation} 
On $Z_1$ vector fields $X^l_i$ and $\tilde{X}_i$ are related by: 
\begin{equation}\begin{split}
    \tilde{X}_i&=-\tilde{x}^i\sum_{k=1}^n \tilde{x}^k X_k^l +\frac{1+|\tilde{x}|^2}{2}X_i^l+\frac{2 \tilde{x}^i}{1+|\tilde{x}|^2} X_{n+1}^l\,,\,\,
    \tilde{X}_{n+1}= \sum_{k=1}^n \tilde{x}^k X_k^l+ \frac{|\tilde{x}|^2-1}{1+|\tilde{x}|^2}X_{n+1}^l\\
    X_i^l &=\frac{2}{1+|\tilde{x}|^2}\left(\tilde{X}_i+\frac{2 \tilde{x}^i}{1+|\tilde{x}|^2}\left(\tilde{X}_{n+1}-\sum_{k=1}^n \tilde{x}^k
        \tilde{X}_k\right)\right)\,,\,\, X_{n+1}^l=\frac{|\tilde{x}|^2-1}{1+|\tilde{x}|^2}\tilde{X}_{n+1}+\frac{2}{1+|\tilde{x}|^2}\sum_{k=1}^n \tilde{x}^k
        \tilde{X}_k,
  \end{split}
\end{equation}
These equalities imply that vector fields $X_i^l\,,\,i=1,\dots, n+1$ are defined 
on $Z$ and $X_i^l(p_0)=0\,,\,i=1,\dots,n$ and $X_{n+1}^l(p_0)=\tilde{X}_{n+1}(p_0)$. 
Commutators of  $\tilde{X}_i$'s are equal:
\begin{equation*}\begin{split}
 [\tilde{X}_i, \tilde{X}_i]&=\frac{4}{1+|\tilde{x}|^2}\left( \tilde{x}^j\tilde{X}_i-\tilde{x}^i\tilde{X}_j\right)\\
[\tilde{X}_i, \tilde{X}_{n+1}]&=\frac{|\tilde{x}|^2-3}{1+|\tilde{x}|^2}\tilde{X}_i-\frac{ 8 \tilde{x}^i}{(1+|\tilde{x}|^2)^2}
\left(\tilde{X}_{n+1}-\sum_{k=1}^n \tilde{x}^k \tilde{X}_k\right)\end{split}
\end{equation*}
To get formulae (\ref{eq:liealg-nawias}) one evaluates   both sides of these equalities at $(x,0,1)$ 
and uses   definition of stereographic coordinates (\ref{eq:stereo}, \ref{eq:stereo-inv}).

It remains to prove formula (\ref{eq:liealg-kotwa}) for the anchor. The anchor at $p\in S^n$ is equal to $T_{(p,0,1)}\tilde{e}_R$ restricted to the
tangent space to the left fiber of $Z$ at $(p,0,1)$. The formula (\ref{eq:inverse}) gives $\tilde{e}_R(p,v,s)=(p_1,0,1)$, where
$p=:\left( \begin{array}{c}u\\ \alpha\end{array} \right)$ and 
$p_1=:\left( \begin{array}{c}u_1\\ \alpha_1\end{array} \right)$ as in (\ref{eq:inverse}). Computing the tangent map one gets:\notka{eq:anchor-Z}
\begin{equation}\label{eq:anchor-Z}
\begin{split}  
  \rho(p,0,1)(\partial_s)&=(p,e_{n+1}-\alpha p)\in T_pS^n\quad,\quad \partial_s:\, ]-\epsilon,\epsilon[\ni t \mapsto (p,0,1+t)\in Z\\
  \rho(p,0,1)(\dot{v})&=(p,(\alpha-1)\dot{v})\in T_pS^n\quad,\quad \dot{v}:\,  ]-\epsilon,\epsilon[\ni t \mapsto (p,\dot{v}t,1)\in Z
\end{split}
\end{equation}
Using the definition of $\tilde{X}_i$ one gets (\ref{eq:liealg-kotwa}).
\dowl

The Lie algebroid structure on $\lie{\tilde{S}}=TS^n$ is as follows. Define sections of $TS^n$ (cf. (\ref{def:basis-sections}))
\begin{equation} \label{def:basis-sections-T} \tilde{Y}_i(p):=(p, P_p(e_i)),\,i=1,\dots, n+1. \end{equation}\notka{prop:liealg-T}
\begin{prop}\label{prop:liealg-T}
  \begin{enumerate}
  \item For any $p\in S^n$:  $T_pS^n=span\{ \tilde{Y}_i(p): i=1,\dots, n+1\}$;
    \item For $p\in  S^n\setminus e_i^\perp\,$ vectors in the set
      $\{\tilde{Y}_1(p),\dots, \tilde{Y}_{n+1}(p)\}\setminus \{\tilde{Y}_i(p)\}$ form a basis in $T_pS^n$;
    \item The Lie algebroid bracket is given by:\notka{eq:bracket-T}
      \begin{equation}\label{eq:bracket-T}
        [ \tilde{Y}_i\,,\,\tilde{Y}_j] =(p^j-\delta_{j, n+1})\tilde{Y}_i -(p^i-\delta_{i, n+1})\tilde{Y}_j \,,\,i.j=1,\dots, n+1\end{equation}
    \item The anchor $\tilde{\rho}:TS^n\rightarrow TS^n$ is equal to:\begin{equation}\label{eq:anchor-T}
        \tilde{\rho}(p,v)=(p, (p^{n+1}-1)v)\,,\,(p,v)\in TS^n, \end{equation}
    \item For any two sections $X, Y$ of $\lie{\tilde{S}}$: $[X,Y](p_0)=0$.
    \end{enumerate}
    \end{prop}

    \noindent{\em Proof:} $1)$ and $2)$ are straightforward;\\ 
    $3)$ Left invariant vector fields on $\tilde{S}_1=T(S^n\setminus\{p_0\})$ corresponding to sections $\tilde{Y}_i$ (in stereographic coordinates) reads: 
\begin{equation} \begin{split} \tilde{Y}_i(x,\dot{x})&=-\frac{1+|x|^2}{1+|\hat{x}|^2} \hat{x}^i\sum_{k=1}^n \hat{x}^k\partial_{\dot{x}^k}+
    \frac{1+|x|^2}{2}\partial_{\dot{x}^i} \,,\, \hat{x}:=x -\frac{2}{1+|x|^2} \dot{x}\,,\,i=1,\dots, n\\
  \tilde{Y}_{n+1}(x,\dot{x})&=\frac{1+|x|^2}{1+|\hat{x}|^2}\sum_{k=1}^n \hat{x}^k\partial_{\dot{x}^k},\end{split}\end{equation} 
and their commutators: 
\begin{equation} [ \tilde{Y}_i\,,\,\tilde{Y}_j]=\frac{2}{1+|\hat{x}|^2} \left( \hat{x}^j \tilde{Y}_i- \hat{x}^i \tilde{Y}_j\right)\,,\,\,
  [ \tilde{Y}_i\,,\,\tilde{Y}_{n+1}]=\frac{-2}{1+|\hat{x}|^2} \left(\tilde{Y}_i+ \hat{x}^i \tilde{Y}_{n+1}\right)\,,\,i,j=1,\dots, n .
\end{equation} 
Evaluating these expressions at $(x,0)$ and using (\ref{eq:stereo}, \ref{eq:stereo-inv}) one obtains brackets  (\ref{eq:bracket-T});\\
$4)$ After a moment of reflections one realizes that the anchor for $\lie{\tilde{S}}$ is given by the second line of (\ref{eq:anchor-Z});\\
$5)$ By  (\ref{eq:bracket-T}): $[ \tilde{Y}_i\,,\,\tilde{Y}_j](p_0)=0$, and by (\ref{eq:anchor-T}): $\tilde{\rho}(p_0,v)=0$. Since sections
$\tilde{Y}_1,\dots,\tilde{Y}_n$ form a basis of sections of $TS^n$ in a neighborhood of $p_0$ the claim follows. 
\dowl   
\begin{re} Sections $\tilde{Y}_i$ in stereographic coordinates are equal:  
  \begin{equation} \tilde{Y}_{n+1}(x)=\sum_{k=1}^n x^k\partial_{x^k}\,,\,
    \tilde{Y}_i(x)=\frac{1+|x|^2}{2}\partial_{x^i}-x^i \tilde{Y}_{n+1}(x)\,,\,i=1,\dots,n
  \end{equation}  
  and the anchor: $\displaystyle \tilde{\rho}(\partial_{x^i})=-\frac{2}{1+|x|^2}\partial_{x^i}$.\\  
  By these formulae and (\ref{eq:bracket-T}) the Lie algebroid bracket reads:  
  \begin{equation} [\partial_{x^i}, \partial_{x^j}]= \frac{4}{(1+|x|^2)^2}\left( x^i \partial_{x^j}- x^j \partial_{x^i}\right)
  \end{equation}  
  These formulae, together with $\tilde{\rho}(X)(p_0)=0$ for any section $X$ and $5)$ of prop.\ref{prop:liealg-T},  
  determine structure of Lie algebroid completely.  
 
  The bracket above define Poisson bracket on $T^*S^n$:   
  \begin{equation} \begin{split} \{F, G\}(\varphi)&=0\,,\,\varphi\in T_{p_0}^*S^n\,,\, F,G\in C^\infty(T^*S^n),\\
\{p_i, p_j\}&:=\frac{4}{(1+|x|^2)^2}\left( x^i p_j- x^j p_i\right)\,,\, \{p_i, f\}=-\frac{2}{1+|x|^2}(\partial_{x^i}f)\,,\,
\{f,g\}=0\,,\,f,g\in C^\infty(S^n), 
\end{split} 
\end{equation}  
where coordinates $(x,p)$ on $T^*(S^n\setminus \{p_0\})$ correspond to stereographic coordinates. Since $\tilde{S}_1$ is isomorphic to a pair groupoid,  
this Poisson bracket is of maximal rank on $T^*(S^n\setminus\{p_0\})$.  
\end{re}  
 
\subsection{Semi-classical limit} 
Our  last task is to justify the name ``quantum space of euclidean lines'' for  the groupoid algebra 
$C_r^*(\tilde{S})$ (at least partially, more complete arguments will
be given in  \cite{PS-future}). It will be shown  that the base map of $T^*\delta_{\tilde{S}}$
(the cotangent lift of morphism $\delta_{\tilde{S}}$, see below) is the action of the euclidean group on the space of oriented euclidean lines
described  in (\ref{eq:action-lines}). In fact we will also show that we get the action (\ref{eq:action-spheres}) on the (Poisson) ``sphere bundle''.  

Recall  that {\em a cotangent lift} of a differentiable relation $h:X\rel Y$ is the relation $T^*h: T^*X\rel T^*Y$ defined as \cite{SZ2}:
$$ T^*h=\{(\psi_y,\varphi_x)\in T_y^*Y\times T_x^*X : \forall (t_y, t_x)\in T_{(y,x)}Gr(h)\,: \,<\psi_y, t_y>=<\varphi_x, t_x>\}$$  
If $X\rra E$ and $Y\rra F$ are differential groupoids then $T^*X\rra (TE)^0$ and $T^*Y\rra (TF)^0$ are their {\em cotangent lifts}. 
If $h:X\rel Y$ is a morphism (of diff. groupoids)  then $T^*h:T^*X\rel T^*Y$ is a morphism as well, so there is the base map $(TF)^0\rightarrow (TE)^0$,
which is a Poisson map \cite{SZ2}.

By (\ref{eq:tildeZ}) the graph of $\delta_{Z}: Z\rel G_B\times Z$ is equal to:\notka{eq:tildeZ-1} 
\begin{equation}\label{eq:tildeZ-1} \delta_{Z}=\{(b \tilde{c}_L(z), z; b z): b\in B, z\in TS^n\times\R_+\},
\end{equation} 
where $z=:(p,v,s)$ and the action $bz=b(p,v,s):=(bp,bv,s)$ (recall that $B=SO(\npj)$), and $\tilde{c}_L: Z\rightarrow C$ was defined 
just before (\ref{eq:tildeZ}) (in particular $\tilde{c}_L(p,0,1)=e_C$, the  neutral element of $C$).
The goal now is to  compute the base map of $T^*\delta_{Z}$. 

By (\ref{eq:tildeZ-1})  the base map of $\delta_{Z}$ is the map  $B\times S^n\ni (b,p)\mapsto bp\in S^n$, so  for the base map of  
   $T^*\delta_{Z}$, denoted by $\beta$, one has: 
   $$ \beta: (T_bB)^0\times (T_{(p,0,1)}(S^n\times\{0\}\times \{1\}))^0\rightarrow (T_{(bp,0,1)}(S^n\times\{0\}\times \{1\}))^0.$$      
     {\em To simplify notation we will write $T_{(p,0,1)}S^n$ instead of $T_{(p,0,1)}(S^n\times\{0\}\times \{1\})$, 
       but keep $(p,0,1)$ to indicate that  it denotes a subspace of $T_{(p,0,1)}Z$}.  
   The submanifold $(TB)^0\subset T^*G$ 
    is a subgroup; by right translation it will be identified with $\gotb^0\rtimes B\subset \gotg^*\rtimes G$; 
i.e.  $(\varphi,b)\in \gotb^0\rtimes B$ represents $\varphi b\in T_b^*G$ so that for $X_b\in T_bG$:  $<\varphi b, X_b>=<\varphi, X_bb^{-1}>$. 
The mapping $\beta$ is defined by:\notka{def:beta} 
    \begin{equation}\label{def:beta}
      \begin{split} \beta(\varphi b, \tilde{\psi}_p)=\tilde{\psi}_{bp}\iff \left[ \, \forall\right.  &  X_b\in T_bG,\, V_p\in T_{(p,0,1)}(TS^n\times \R_+),\,
          V_{bp}\in T_{(bp,0,1)}(TS^n\times \R_+) : \\
        &  \left. (X_b, V_p; V_{bp})\in T\delta_{Z} \,\,\Rightarrow \,\,<\tilde{\psi}_{bp},V_{bp}>= <\varphi b , X_b> +
        <\tilde{\psi}_p,V_{p}>\,\right] \end{split}
      \end{equation}
By (\ref{eq:tildeZ-1}), condition $(X_b, V_p; V_{bp})\in T\delta_{Z}$ means there are curves
$z(t)=(p(t), v(t), s(t))\in Z\,,\,z(0)=(p, 0, 1)$ representing $V_p$,
$b(t)\in B\,,\,b(0)=b$ with $b(t)\tilde{c}_L(z(t))$ representing $X_b$, and  the vector  $V_{bp}$ is represented by
$b(t)z(t)=(b(t) p(t), b(t) v(t), s(t))$.

Let us begin by the term $<\tilde{\psi}_{bp},V_{bp}>$. The vector $V_{bp}$ is represented by $b(t)z(t)$ i.e.
$V_{bp}=\dot{b}z(0)+ bV_p$, where $\dot{b}z(0)$ is tangent to the mapping $B\ni b\mapsto bz(0)\in TS^n\times \R_+$
applied to $\dot{b}\in T_bB$. Since $z(0)=(p,0,1)$, $bz(0)=(bp,0,1)\in S^n\times\{0\}\times \{1\}$ and  $\dot{b}z(0)\in T_{(bp,0,1)}S^n$,
but $\tilde{\psi}_{bp}\in (T_{(bp,0,1)}S^n)^0$ so that $\dsp <\tilde{\psi}_{bp}, V_{bp}>=<\tilde{\psi}_{bp}, b V_p>.$

Let us consider the term $<\varphi b , X_b>$. For $p\in S^n$ let us denote  $F_p:=T_{(p,0,1)}\tilde{c}_L: T_{(p,0,1)}Z\rightarrow \gotc$. 
We will show the equality:\label{eq:var-fp}
\begin{equation}\label{var-fp}
  <\varphi b , X_b>=<\varphi , X_bb^{-1}>=<\varphi,F_{bp}(b V_p)>=<F_{bp}^*(\varphi), b V_p>.
  \end{equation}
The vector
$X_bb^{-1}\in \gotg$ is represented by the curve $b(t)\tilde{c}_L(z(t))b^{-1}$; since
$\tilde{c}_L(z(0))=e_C$  and $\varphi\in \gotb^{0}$:
$$<\varphi , X_bb^{-1}>=<\varphi , bT_{z(0)}\tilde{c}_L(V_p)b^{-1}>=<\varphi ,P_{\gotc}[\ad(b)(T_{z(0)}\tilde{c}_L(V_p))]>=
<\varphi ,P_{\gotc}[\ad(b)(F_p(V_p))]>,$$
where $\ad()$ is the adjoint representation of $G$ on $\gotg$ and $P_{\gotc}$ is a projection defined by the (vector space)
decomposition $\gotg=\gotb\oplus\gotc$. This way (\ref{var-fp}) will follow from:
$\dsp P_{\gotc}[\ad(b)(F_p(V_p))]=F_{bp}(b V_p)$.

For $b\in B$ and $z\in Z$ it is clear that 
$\tilde{c}_L(bz)=c_L(b\tilde{c}_L(z))$ so by applying the tangent mapping:
$\dsp T_{bz}\tilde{c}_L(b V)=T_{b\tilde{c}_L(z)} c_L(b T_{z}\tilde{c}_L(V))\,,\,V\in T_{z}Z$.
For $z:=z(0)=(p,0,1)$ we have $bz=(bp,0,1)$ (and $\tilde{c}_L(z)=e_C$). Thus:
$$F_{bp}(b V_p)=T_{(bp,0,1)}\tilde{c}_L(b V_p)=T_{b} c_L(b F_p(V_p)),$$
but since $c_L(g)=c_L(g b^{-1})$, we have, for any $X_b\in T_bG$ the equality: $T_bc_L(X_b)=P_{\gotc}(X_b b^{-1})$, so
$$T_{b} c_L(b F_p(V_p))=P_{\gotc}(b  F_p(V_p)b^{-1})= P_{\gotc}[\ad(b)(F_p(V_p))].$$
This way  (\ref{var-fp}) is proven and (\ref{def:beta}) can be rewritten as:
      \begin{equation*}
      \beta(\varphi b, \tilde{\psi}_p)=\tilde{\psi}_{bp}\iff \left[\, \forall   V_p\in T_{(p,0,1)}(TS^n\times \R_+):\,
        <\tilde{\psi}_{bp}, b V_{p}>=<F_{bp}^*(\varphi), b V_p> + <\tilde{\psi}_p,V_{p}>\,\right] 
    \end{equation*}  
    and finally we get the formula for $\beta$:\notka{eq:Eaction}  
    \begin{equation}
      \label{eq:Eaction}
      \beta(\varphi b, \tilde{\psi}_p)=b\tilde{\psi}_p+ F_{bp}^*(\varphi)\,,\,\,b\in B\,,\,\varphi\in \gotb^0\subset \gotg^*\,,\,
      \tilde{\psi}_p\in (T_{(p,0,1)}S^n)^0\subset T^*_{(p,0,1)}Z
      \end{equation}  
      Since the morphism $\delta_{Z}$ satisfies (\ref{eq:coact-tildez}):
      $(\delta_B\times id)\delta_{Z}=(id\times \delta_{Z})\delta_{Z}$,
      by the general theory of cotangent lifts \cite{SZ2} this formula defines Poisson action of a Poisson-Lie group   $\gotb^0\rtimes B$ 
 on the Poisson manifold $(T(S^n\times\{0\}\times\{1\}))^0\subset T^*Z$ 
    (which is the dual bundle to the Lie algebroid of $Z$). 

Let us now show, that after some identifications, this is the formula ({\ref{eq:action-spheres}).  The bundle $(T(S^n\times\{0\}\times\{1\}))^0$
can be identified with $T^*S^n\oplus \bunR$ (which, in turn will be identified with $TS^n\oplus \bunR$) by the following considerations.   
Clearly there is a decomposition: $T_{(p,v,s)}(TS^n\times \R_+)=T_{(p,v)}TS^n\oplus T_s\R=T_{(p,v)}TS^n\oplus \R$, 
moreover for $v=0$ we have the decomposition $T_{(p,0)}(TS^n)=T_pS^n\oplus T_pS^n$; elements of the first summand are vectors from $T_{(p,0)}S^n$
(i.e. $(p,v)$ is represented by a curve $t\mapsto (p(t),0)\in TS^n$  with $p(0)=p\,,\,\dot{p}(0)=v$), and   the second  summand is the space
tangent to the fiber (i.e. $(p,v)$ is represented by a curve $t\mapsto (p,tv)\in TS^n$). So we can write:
$$T_{(p,0,s)}(TS^n\times \R_+)=T_pS^n\oplus T_pS^n\oplus \R$$
and $\psi_p+\rho\in T_p^*S^n\oplus \R$ as element of $(T_{(p,0,1)}(S^n\times\{0\}\times\{1\}))^0$ acts as
$$<\psi_p+\rho, w + v + z>:=<\psi_p,v>+ \rho z\,,\, w + v+ z\in T_pS^n\oplus T_pS^n\oplus \R=T_{(p,0,1)}(TS^n\times \R_+)$$
Since $\tilde{c}_L(p,0,1)=e_C$, for the mapping $F_p$ we have $F_p(w\oplus 0\oplus 0)=0$, so $F_p$ could and {\em will be considered from now on}
as a mapping $T_pS^n\oplus\R\rightarrow \gotc$; consequently $F_p^*:\gotc^*=(\gotb)^0\rightarrow T^*_pS^n\oplus \R$.

Now $F_p$ will be computed explicitly. Let $(p,v)\in TS^n$ and
$\displaystyle p=\left(\begin{array}{c} u\\\alpha\end{array}\right),\, v=\left(\begin{array}{c} \mu\\ r \end{array}\right)$,
with $\mu^t u + r\alpha=0$. Then $\tilde{c}_L(p,v,s)=(\tilde{s},\tilde{y})\in C$, where
     $(\Lambda,u,w,\alpha; s, y)=(\tilde{s}, \tilde{y}) (\tilde{\Lambda}, \tilde{u}, \tilde{w}, \tilde{\alpha})$ (this is decomposition (\ref{eq:iwasawa})),
     and  $\mu=\Lambda y\,,\,r=w^t y$. The solution (see \cite{PS-kappa2})  gives:\notka{eq:tildecL} 
     \begin{equation}\label{eq:tildecL}
       \tilde{c}_L\left(\left(\begin{array}{c} u\\\alpha\end{array}\right),\left(\begin{array}{c} \mu\\ r \end{array}\right), s\right)=
       \left(\tilde{s}:=s -r +(\alpha-1)\frac{s^2-1-|v|^2}{2s},\,\tilde{y}:=\mu- \frac{s^2-1-|v|^2}{2s} u\right)\end{equation} 
     Let $(p,\dot{v},\dot{s})\in T_pS^n\oplus \R\,,\,\dot{v}:=\left(\begin{array}{c} \dot{\mu}\\ \dot{r} \end{array}\right)$;
     as a vector tangent to the fiber $T_{(p,0)}\oplus \R$ at $(p,0,1)$ it is represented by the curve $]-\epsilon,\epsilon[\ni t \mapsto (p, t\dot{v}, 1+\dot{s}t)=:(p, v(t), s(t)) \in TS^n\times\R_+$, therefore  $F_p(p,\dot{v},\dot{s})$ is represented by the curve:
     $$]-\epsilon,\epsilon[\ni t \mapsto \left(\tilde{s}(t)=s(t) -r(t) +(\alpha-1)\frac{s(t)^2-1-|v(t)|^2}{2s(t)}\,,\,
       \tilde{y}(t):=\mu(t)-\frac{s(t)^2-1-|v(t)|^2}{2s(t)} u \right)$$
     By differentiation one obtains:\notka{eq:Fp}
     \begin{equation}\label{eq:Fp}
       F_p(p,\dot{v},\dot{s})=(-\dot{r}+\alpha \dot{s})\partial_s+\sum_{k=1}^n (\dot{\mu}^k-\dot{s}u^k)\partial_{y^k}
       \end{equation}
    
       To identify $\gotb^0\rtimes B$ with $V\rtimes SO(\eta)$ for an Euclidean vector space $(V,\eta)$ we present
       the Iwasawa decomposition (\ref{eq:iwasawa}) in a more geometric way
       (cf.  \cite{PS-poisson}), this way a geometric meaning of various
       identifications can be clearly seen.
       Let  $(V,\tilde{\eta})$ be a  vector Minkowski space of dimension $n+2\,,\,n>1$
       (signature of  $\tilde{\eta}$ is $(+,-,\dots, -)$) and  $G:=SO_0(\tilde{\eta})$ and consider the semidirect product
       $\gotg^*\rtimes G$:\notka{semi-direct}
\begin{equation}
\label{semi-direct}
(\varphi,g)(\psi, h):=(\varphi+\kad(g) \psi, g h),
\end{equation} 
where $\gotg^*$ is the dual of the Lie algebra of $G$ and $\kad(g)$ is the coadjoint representation. The Lie algebra of $G$ is 
$\gotg:=so(\tilde{\eta})=span\{M_{xy}:x,y\in V\}$, where
$M_{xy}:=x\mt\tilde{\eta}(y)-y\mt\tilde{\eta}(x)$ and $\tilde{\eta}:V\rightarrow V^*$ is the isomorphism defined be the form $\tilde{\eta}$.
$\gotg$ is equipped with the non-degenerate bilinear form
$k(M_{xy}, M_{zt}):=\tilde{\eta}(x,t)\tilde{\eta}(y,z)-\tilde{\eta}(x,z)\tilde{\eta}(y,t)$ (the Killing form up to a constant) and
by the associated isomorphism $k:\gotg\rightarrow \gotg^*$ there is the non-degenerate form $\tilde{k}$ on $\gotg^*$ which satisfies:
$\tilde{k}(\kad(g) k(X),k(Y))= k(\ad(g) X, Y)\,\,,\,X,Y\in so(\tilde{\eta})$; it follows $\kad(g)\in O(\tilde{k})$.
Choose two (orthogonal) vectors $\tgr, \sgr\in V$ with
$\tilde{\eta}(\tgr, \tgr)=1=-\tilde{\eta}(\sgr, \sgr)\,,\,\tilde{\eta}(\tgr, \sgr)=0$, and let $\fgr:=\tgr-\sgr$. The orthogonal complement $\tgr^\perp$
with the form $\eta:=-\tilde{\eta}$ is $n+1$ dimensional euclidean space  as well as $(\gotb^0,\tilde{k}|_{\gotb^0}) $ and the map
$M_{\tgr}: \tgr^\perp\ni w\mapsto k(M_{\tgr w})\in \gotb^0$ is an isometry.  The subspace $\gotb:=\{M_{xy}: x, y \in \tgr^\perp\}$ is a Lie subalgebra and
its Lie subgroup $B$ - the stabilizer of $\tgr$ - is $SO(\eta)$. The map $M_{\tgr}\times id$ is an
isomorphism of $\tgr^\perp\rtimes B$ with the subgroup $ \gotb^0\rtimes B\subset \gotg^0\rtimes G$. 
The subspace $\gotc:=span\{M_{\fgr x}: x\in  \tgr^\perp\}$ is a Lie subalgebra and the set
$C:=\{\exp(\log s M_{\sgr\tgr})\exp (M_{\fgr y}): s\in\R_+, y\in <\sgr, \tgr>^\perp\}$ is a Lie group which is a semidirect product of
$C_1:=\{\exp(\log s M_{\sgr\tgr})\simeq \R_+$ and the commutative group $\{\exp (M_{\fgr y}): y\in <\sgr, \tgr>^\perp\}$.
The Iwasawa decomposition is  $SO(\tilde{\eta})=BC$.
The subgroupoid $Z=TS^n\times C_1\subset \tgr^\perp\times  \tgr^\perp\times C_1$ with the projection 
$\pi_Z: G\ni b \exp(\log s M_{\sgr\tgr}) \exp (M_{\fgr y})\mapsto (p:=b\sgr, v:=by, \exp(\log s M_{\sgr\tgr}))$.  
All our previous computations were  made after the  choice of an orthonormal basis in $V$ 
with  $e_0=\tgr$, $p_0=e_{n+1}=\sgr$ (i.e. $S^n\subset \R^{\npji}=\tgr^\perp$) and after the identification of $C_1$ with $\R_+$. 
Finally we identify $TS^n$ with $T^*S^n$ by $TS^n\ni (p,v)\mapsto (p,\eta(v))\in T^*S^n$.  
With all these identifications, formulae 
(\ref{eq:Eaction}) and (\ref{eq:action-spheres}) are equivalent if the following equality holds:  \notka{eq:equiv}
\begin{equation}\label{eq:equiv}\begin{split}
    \forall \, & z\in \tgr^\perp, \, b\in B,\, (p,\dot{v})\in TS^n\subset \tgr^\perp\times \tgr^\perp,\,\dot{s}\in \R:\\
    & \beta(k(M_{\tgr z})b; p, \eta(\dot{v}); \dot{s} k(M_{\tgr \sgr}))=(bp,\eta(b\dot{v}+z -\eta(z, bp) bp); (\dot{s}+\eta(z,bp))k(M_{\tgr \sgr})))
    \end{split}
  \end{equation} 

  With $(s,y)=\exp(\log s M_{\sgr\tgr})\exp (M_{\fgr y})$ it is easy to check that $\partial_s=M_{\sgr \tgr}$ and  
$\partial_{y^k}=M_{\fgr e_k}\,,\,k=1,\dots, n$. 
  Writing $p=u+\alpha\sgr,\,u\in <\sgr, \tgr>^\perp$ and $\dot{v}:=\dot{\mu}+\dot{r}\sgr,\, \eta(p,\dot{v})=0$ (so that $(p,\dot{v})\in TS^n$) the
  formula (\ref{eq:Fp}) reads:
  $$F_p(\dot{v}, -\dot{s}\partial_s)=F_p(\dot{v}+\dot{s} M_{\tgr \sgr})=-(\dot{r}+\alpha\dot{s} )M_{\sgr \tgr}+M_{\fgr,\dot{\mu}+\dot{s} u}=
  (\dot{r}+\alpha\dot{s} )M_{\tgr \sgr}+M_{\fgr,\dot{\mu}+\dot{s} u}$$
  and for $F_p^*$ (after identification $\gotb^0\simeq \gotc^*$), note also that $\eta(\dot{v},p)=0$: 
  \begin{equation*}
    \begin{split}
      <F_p^*(k(M_{\tgr z})), \dot{v}+\dot{s}M_{\tgr \sgr}>& =<(k(M_{\tgr z}), F_p( \dot{v}+\dot{s}M_{\tgr \sgr})>=
      k(M_{\tgr z}, (\dot{r}+\alpha\dot{s} )M_{\tgr\sgr}+M_{\fgr,\dot{\mu}+\dot{s} u})=\\
      &=\eta(z,\dot{v})+\dot{s}\eta(z,p)=\eta(z-\eta(z,p)p,\dot{v})+\eta(z,p) <k(M_{\tgr \sgr}),\dot{s} M_{\tgr \sgr}>=\\
      &=<\eta(z-\eta(z,p)p)+\eta(z,p) k(M_{\tgr \sgr}), \dot{v}+\dot{s} M_{\tgr \sgr}>
    \end{split}
  \end{equation*} 
  Since $\eta(z-\eta(z,p)p,p)=0$ we get 
  $$F_p^*(k(M_{\tgr z}))=(p,\eta(z-\eta(z,p)p);\eta(z,p) k(M_{\tgr \sgr}))\in T_p^*S^n\oplus \gotc_1^*.$$      
  This way for the formula (\ref{eq:Eaction}) written for $\tilde{\psi}_p=(p,\eta(\dot{v});\dot{s}k(M_{\tgr \sgr})\in T_p^*S^n\oplus \gotc_1^*$: 
  \begin{equation*}
    \begin{split}
      \beta(k(M_{\tgr z})b; p, \eta(\dot{v}); \dot{s} k(M_{\tgr \sgr}))& =b(p, \eta(\dot{v}); \dot{s} k(M_{\tgr \sgr}))+F_{bp}^*(k(M_{\tgr z}))=\\ 
       &= (bp,\eta(b\dot{v});\dot{s} k(M_{\tgr \sgr}))+(bp,\eta(z-\eta(z,bp)bp);\eta(z,bp) k(M_{\tgr \sgr}))=\\
       &= (bp, \eta(b\dot{v}+z-\eta(z,bp)bp); (\dot{s}+\eta(z,bp)) k(M_{\tgr \sgr})) 
      \end{split}
    \end{equation*}
    and this is the equality (\ref{eq:equiv}).

    If we look at calculations and formulae starting from (\ref{eq:tildeZ-1}) we observe that:\\
    $(i)\,$  the formula (\ref{eq:tildeZ-1}) defines (by restriction)  the morphism
    $\delta_{\tilde{S}}:\tilde{S}\rel G_B\times \tilde{S}$ and $(\delta_B\times id)\delta_{\tilde{S}}=(id\times \delta_{\tilde{S}})\delta_{\tilde{S}}$;\\
    $(ii)\,$ thus the base map of $T^*\delta_{\tilde{S}}$ defines the action of a Poisson-Lie group $\gotb^0\rtimes B$ 
    on the Poisson manifold $(T(S^n\times\{0\}))^0\subset T^*\tilde{S}$ which, as above, can be identified with $T^*S^n$;\\
    $(iii)\,$ to get the formulae for $T^*\delta_{\tilde{S}}$ it is enough to put $\dot{s}=0$ in formulae for  $T^*\delta_{Z}$;\\
    $(iv)\,$ this way, after  identifications as above, one obtains the action  (\ref{eq:action-lines}).

  
\end{document}